\documentclass[11pt,twoside,a4paper]{article}
\usepackage{amsmath}
\usepackage{amsfonts}
\usepackage{amssymb}
\usepackage{latexsym}
\usepackage{graphics}
\usepackage{theorem}
\usepackage{epsf}
\usepackage{epsfig}
\usepackage[english]{babel}
\usepackage{afterpage}
\usepackage[all]{xy}
\usepackage{hhline}
\usepackage{epic,eepic}

\newtheorem{theorem}{Theorem}

\newtheorem{remark}{Remark}
\newtheorem{algo}{Algorithm}
\newtheorem{lemma}{Lemma}

\newenvironment{@abssec}[1]{%
     \if@twocolumn
       \section*{#1}%
     \else
       \vspace{.05in}\footnotesize
       \parindent .2in
         {\upshape\bfseries #1. }\ignorespaces
     \fi}
     {\if@twocolumn\else\par\vspace{.1in}\fi}
\renewenvironment{abstract}{\begin{@abssec}{\abstractname}}{\end{@abssec}}
\renewcommand\abstractname{Abstract}
\newenvironment{keywords}{\begin{@abssec}{\keywordsname}}{\end{@abssec}}
\newcommand\keywordsname{Key words}
\newenvironment{AMS}{\begin{@abssec}{\AMSname}}{\end{@abssec}}
\newcommand\AMSname{AMS subject classifications}
\newenvironment{mythanks}{\begin{@abssec}{\thanksname}}{\end{@abssec}}
\newcommand\thanksname{Acknowledgements}

\def\D{\mathop {\cal D}}
\def\X{\mathop {\cal X}}

\title{Reduced-Basis approach for homogenization beyond the periodic setting}

\date{February 22, 2007}

\author{S\'ebastien \textsc{Boyaval} \footnote{\textbf{CERMICS}, Ecole Nationale des Ponts et Chauss\'ees, 6 \& 8 avenue Blaise Pascal, Cit\'e Descartes, 77455 Marne-la-Vall\'ee Cedex 2, France and \textbf{MICMAC Project}, INRIA, Domaine de Voluceau, BP. 105 - Rocquencourt, 78153 Le Chesnay Cedex (boyaval@cermics.enpc.fr)}}
\begin{document}

\maketitle

\begin{abstract}
We consider the computation of averaged coefficients for the homogenization of elliptic partial differential equations. In this problem, like in many multiscale problems, a large number of similar computations parametrized by the macroscopic scale is required at the microscopic scale. This is a framework very much adapted to model order reduction attempts. 

The purpose of this work is to show how the reduced-basis approach allows to speed up the computation of a large number of cell problems without any loss of precision. The essential components of this reduced-basis approach are the {\it a posteriori} error estimation, which provides sharp error bounds for the outputs of interest, and an approximation process divided into offline and online stages, which decouples the generation of the approximation space and its use for Galerkin projections.
\end{abstract}



\tableofcontents

\section{Introduction}

In this work, we study the numerical homogenization of \emph{linear scalar elliptic} partial differential equations (PDEs) such as those encountered in the problems of thermal diffusion and electrical conduction. \emph{Oscillating test functions}, also termed \emph{correctors}, are computed through a \emph{reduced-basis} (RB) approach for parametrized \emph{cell problems} supplied with periodic boundary conditions. Numerical results have been obtained with some prototypical parametrizations of the oscillating coefficients and are shown in a two-dimensional case with one single varying rectangular inclusion inside rectangular cells. The method applies to all numerical homogenization strategies that require to solve a large number of parametrized cell problems.

In periodic homogenization, only \emph{one} cell problem has to be solved in order to completely determine the homogenized coefficient(s) to be used in the homogenized (macroscopic) equation. In sharp contrast, non-periodic homogenization requires the solution of several cell problems (in fact, theoretically, an \emph{infinite} number of them, and in practice, \emph{a large number}). A homogenized coefficient is then approximated by some average over a large number of microscopic cells. Consequently, as opposed to the periodic case where the computation is light and exact, the non-periodic case asks for a computationally demanding and approximate-in-nature task. This is why the design of a fast and accurate numerical homogenization method is considered as an important issue for the treatment of non-periodic heterogeneous structures. The RB approach seems very well adapted to this framework. 

The article is organized as follows. In section 2, we give a detailed presentation of the setting of the problem. For the sake of consistency and the convenience of the reader, we also briefly outline the main relevant issues in homogenization and RB theories. In section 3, the RB approach for a parametrized cell problem is introduced and we notably derive {\it a posteriori} error bounds related to the convergence of the RB method in the homogenization context. Numerical results for the prototypical example of rectangular cells with one single rectanguler inclusion are presented in section 4. Possible extensions of our work are discussed in the final section.

\section{Setting of the problem, elements of homogenization theory and RB approach}

\subsection{Formulation of the problem}

The mathematical problem under consideration throughout this article reads as follows. We are interested in the behaviour of a sequence of scalar functions $u^{\epsilon}$ that satisfy
\begin{equation} \label{pbeps}
- {\rm div}( \bar{\bar{A}}^{\epsilon}(x) \nabla u^{\epsilon}(x) ) = f(x), \forall x \in \Omega
\end{equation}
in a bounded open set $\Omega \subset \mathbb{R}^{n}$, for a sequence of scalars $\epsilon > 0$. Of interest is the asymptotic limit of the sequence $u^{\epsilon}$ when $\epsilon \rightarrow 0$, along with approximations for $u^{\epsilon}$ when $\epsilon$ is small. 

For the sake of simplicity, the scalar source term $f$ is chosen in $L^{2}(\Omega)$ and we supply equation (\ref{pbeps}) with the following boundary conditions on the smooth (say ${\cal C^{1}}$-Lipschitz) boundary $\partial \Omega = \Gamma_{D} \bigcup \Gamma_{N}$ of $\Omega$,
\begin{equation} \label{BC}
(BC) \left\{ \begin{array}{l}
u^{\epsilon} \mid_{\Gamma_{D}} = 0 \\
\bar{\bar{A}}^{\epsilon} \nabla u^{\epsilon}\cdot\bar{n} \mid_{\Gamma_{N}} = 1\ .
\end{array} \right.
\end{equation}
But as a matter of fact, it is well known that the homogenization results are local in nature and do not depend on the boundary conditions, except for what regards error estimations close to the boundary. Nor do the homogenization results depend on the source term $f$. Hence the generality of the assumptions (BC) and $f \in L^{2}(\Omega)$, chosen here to give a precise mathematical frame to the numerical experiments. 

To fix ideas, the unknown $u^{\epsilon}$ could be thought of, either as a temperature or as an electric field in a macroscopic domain $\Omega$. The tensorial coefficients for $\bar{\bar{A}}^{\epsilon}(x)$ would respectively be thought of, either as temperature diffusivities or as electric conductivities.

Next, let us define, for any $\epsilon > 0$, the family $\bar{\bar{A}}^{\epsilon} \in L^{\infty}(\Omega,{\cal M}_{\alpha_{A},\gamma_{A}})$ of functions from $\Omega$ to the set ${\cal M}_{\alpha_{A},\gamma_{A}}$ of uniformly positive definite $n \times n$ matrices (second order tensors) with uniformly positive definite inverses, that is, matrices $\bar{\bar{A}}^{\epsilon}$ satisfying, for all $x \in \Omega$,
\begin{eqnarray} \label{coercivity-continuity}
 0 & < \alpha_{A} \mid u \mid^{2} \leq & \bar{\bar{A}}^{\epsilon}(x) u \cdot u
, \forall u \in \mathbb{R}^{n} \\
 0 & < \gamma_{A} \mid u \mid^{2} \leq & {{\bar{\bar{A}}^{\epsilon}}(x)}^{-1} u \cdot u
, \forall u \in \mathbb{R}^{n}\ .
\end{eqnarray}

Under such conditions, equations \eqref{pbeps}-\eqref{BC} are well posed in the sense of Hadamard. For every $\epsilon > 0$, there exists a unique solution $u^{\epsilon}$ in  $H^{1}_{\Gamma_{D}}(\Omega) = \left\{ u \in H^{1}(\Omega), u \mid_{\Gamma_{D}} = 0 \right\}\ $ that continuously depends on $f$,
\begin{equation}
\| u^{\epsilon} \|_{H^{1}(\Omega)} \leq C(\Omega) \|f\|_{L^{2}(\Omega)}\ ,
\end{equation}
with some constant $C(\Omega)$ that is only function of $\Omega$. Moreover, the sequence of solutions $u^{\epsilon}$ is bounded in $H^{1}_{\Gamma_{D}}(\Omega)$, so that some subsequence $\epsilon'$ weakly converges to a limit $u^\star \in H^{1}_{\Gamma_{D}}(\Omega)$ when $\epsilon' \rightarrow 0$. We are specifically interested in estimating the behaviour of this weakly-convergent subsequence.

In a typical frame for the  homogenization theory, the coefficients $\bar{\bar{A}}^{\epsilon}$ are assumed to oscillate very rapidly on account of numerous small heterogeneities in the domain \nolinebreak[4] $\Omega$. For example, $\epsilon$ typically denotes the ratio of the mean period for microscopic fast oscillations of $\bar{\bar{A}}^{\epsilon}$ divided by the mean period for macroscopic slow oscillations of $\bar{\bar{A}}^{\epsilon}$ in $\Omega$. Moreover, it is usually assumed that macroscopic (macro) and microscopic (micro) scales ``separate'' when $\epsilon$ is sufficiently small, which allows for the oscillating coefficients to be locally \emph{homogenized} in the limit $\epsilon \rightarrow 0$.

\subsection{General context for homogenization}
\label{homogenization_context}

As announced above, this section \ref{homogenization_context} includes some basics of homogenization theory for linear scalar elliptic PDEs. The purpose of this summary is only to collect some elementary results for convenience. Readers familiar with the homogenization theory may then like to skip this section and proceed to section \ref{RB_context}, which introduces the RB theory.

\subsubsection{Abstract homogenization results}

The following abstract homogenization result is the basis for many studies that aim at computing a numerical approximation for $u^{\epsilon}$ when $\epsilon$ is small \cite{Hou-97,Matache-99,Allaire-04,Gloria-06-1,E-07}. It shows that, in the limit $\epsilon \rightarrow 0$, the small oscillating scale ``disappears" from the macroscopic point of view ; that is, the microscopic and macroscopic behaviours asymptotically ``separate". This implies that the limit problem is easier to solve than equation \eqref{pbeps} for some small $\epsilon$, since the former does not require to resolve microscopic details. Moreover, a tractable approximation of $u^\epsilon$ when $\epsilon$ is small enough can be computed from the asymptotic limit when $\epsilon \rightarrow 0$.

More precisely, $u^\star$ can be obtained as the solution to the \emph{H-limit} equation for \eqref{pbeps} (see equation \eqref{pbstar} below). It is then an $L^2$-approximation for $u^{\epsilon}$ when $\epsilon$ is small, as the asymptotic $L^2$-limit of $u^{\epsilon}$ when $\epsilon \rightarrow 0$. Moreover, an improved $H^1$-approximation for $u^{\epsilon}$ when $\epsilon$ is small can also be computed with $u^\star$ after ``correction" of the gradient $\nabla u^\star$. 

The homogenization of the sequence of equations (\ref{pbeps}) is the mathematical process which allows to define the H-limit equation and the $H^1$ approximation for $u^{\epsilon}$. It is performed using the following abstract objects \cite{Jikov-94}: 
\begin{description}
\item[$\bullet$] a sequence of $n$ \emph{oscillating test functions} $z^{\epsilon}_{i} \in H^{1}(\Omega)$ such that, for every direction $(e_i)_{1\leq i\leq n}$ of the ambient physical space $\mathbb{R}^n$, we have $ z^{\epsilon}_{i} {\rightharpoonup} x_{i} \  \text{in}\ {H^{1}(\Omega)}$
and
$$ - {\rm div}( \bar{\bar{A}}^{\epsilon} \nabla z^{\epsilon}_{i} ) = - {\rm div}( \bar{\bar{A}}^{\star} e_{i}) \  \text{in}\ H^{-1}(\Omega)\ ,$$
\item[$\bullet$] a \emph{homogenized tensor} $\bar{\bar{A}}^\star$ defined by
\begin{equation} \label{homogenized_tensor}
\bar{\bar{A}}^{\epsilon} \nabla z^{\epsilon}_{i} \rightharpoonup \bar{\bar{A}}^{\star} e_{i}\  \text{in}\ [L^2(\Omega)]^{n}\ ,
\end{equation}
\item[$\bullet$] a subsequence $u^{\epsilon'}$ of solutions for (\ref{pbeps}) that satisfies
\begin{equation} 
\left\{ \begin{array}{l} 
u^{\epsilon'} \rightharpoonup u^{\star} \ \text{in}\ H^{1}_{\Gamma_{D}}(\Omega) \\
\bar{\bar{A}}^{\epsilon'} \nabla u^{\epsilon'} \rightharpoonup \bar{\bar{A}}^{\star} \nabla u^{\star} \ \text{in}\ [L^{2}(\Omega)]^{n}
\end{array} \right.
\end{equation}
where $u^{\star}$ is solution for the H-limit or \emph{homogenized} equation
\begin{equation} \label{pbstar}
- {\rm div}( \bar{\bar{A}}^{\star}(x) \nabla u^{\star}(x) ) = f(x),\ \forall x \in \Omega\ ,
\end{equation}
supplied with the boundary conditions (BC),
\item[$\bullet$] and an asymptotic approximation for a subsequence $\epsilon'$ of $\epsilon$ that satisfies
\begin{eqnarray} 
\label{L2approx}
\left\| u^{\epsilon'} - u^{\star} \right\|_{L^{2}(\Omega)} & \stackrel{\epsilon' \rightarrow 0}{\longrightarrow} & 0 \\
\label{H1approx}
\left\| \nabla u^{\epsilon'} - \sum_{i=1}^{n} z^{\epsilon'}_{i} \partial_{i} u^{\star} \right\|_{[L^{1}_{loc}(\Omega)]^{n}} & \stackrel{\epsilon' \rightarrow 0}{\longrightarrow} & 0\ ,
\end{eqnarray}
where $\partial_{i} u^{\star}$ are the components of $\nabla u^\star$ in each direction $e_i$.
\end{description}

Note that the latter convergence result (\ref{H1approx}) for $\nabla u^{\epsilon'}$ also holds in $[L^{2}_{loc}(\Omega)]^{n}$ if $u^\star \in W^{1,\infty}(\Omega)$. So, if $u^\star \in H^2(\Omega)$, the \emph{corrector result} states that $u^{\epsilon}$ can be approximated with the following formula,
\begin{equation} \label{corrector}
u^\epsilon = u^\star + \mathop{\sum}_{i=1}^{n} (z^{\epsilon}_{i}-x_i) \partial_i u^\star + r_\epsilon\ ,
\end{equation}
where the remainder term $r_\epsilon$ converges strongly to zero in $W^{1,1}_{loc}(\Omega)$.

In a nutshell, the homogenization of the sequence of equations \eqref{pbeps} has allowed to derive an abstract homogenized problem, \eqref{homogenized_tensor}-\eqref{pbstar}, the solution $u^{\star}$ of which can be corrected with \eqref{H1approx} into an $H^1$ approximation of $u^{\epsilon}$ in the limit $\epsilon \rightarrow 0$. 

But we lack an explicit expression for the homogenized tensor $\bar{\bar{A}}^\star$ to get an explicit asymptotic limit $u^\star$. That is why, though it is not required by the previous abstract theory, the scale separation in the behaviour of the oscillating coefficients $\bar{\bar{A}}^{\epsilon}$ is often assumed to be explicitly encoded, using some specific postulated form for $\bar{\bar{A}}^{\epsilon}$. This allows to derive an explicit expression of the homogenized problem, and even an error estimate in terms of $\epsilon$ for the \emph{correction error} $r_\epsilon$ in \eqref{corrector}, which allows to quantify the homogenization approximation error. 

\subsubsection{The explicit two-scale homogenization}

To get explicit expressions for the homogenized problem, some particular dependence of the family $\bar{\bar{A}}^{\epsilon}$ on the space variable $x$ is often assumed, like in two-scale homogenization for instance. Namely, on account of the scale separation assumption and the local dependence of the homogenization process, one of the most common assumption is the \emph{local periodicity} for $\bar{\bar{A}}^{\epsilon}$, which can be made precise as follows. 

It is assumed that tensors $\bar{\bar{A}}^{\epsilon}$ are traces of functions of two coupled variables on the set locally defined by a fast microscopic variable ${\epsilon}^{-1}x$ linearly coupled with the slow macroscopic variable $x$ in $\Omega$:
\begin{equation} \label{twoscale}
\bar{\bar{A}}^{\epsilon}(x) = \bar{\bar{A}}\left(x,\frac{x}{\epsilon}\right)\ ,
\end{equation}
where, for any $x \in \Omega$, the function $\bar{\bar{A}}(x,\cdot)$
$$\bar{\bar{A}}(x,\cdot):y\in \mathbb{R}^{n}\rightarrow \bar{\bar{A}}(x,y)\in \mathbb{R}^{n \times n}$$ 
is 1-periodic in each of the $n$ directions $(e_{i})_{1 \leq i \leq n}$, which makes the local oscillations completely determined when $\epsilon \rightarrow 0$. The domain $Y = [0,1]^n$ of the periodic pattern is called the \emph{cell} and is identified with the $n$-dimensional torus. $\bar{\bar{A}}(x,\cdot)$ is said to be \emph{Y-periodic}. Note that the properties of the tensors $\bar{\bar{A}}^{\epsilon}$ imply $\bar{\bar{A}} \in L^\infty(\Omega,L^\infty(Y,{\cal M}_{\alpha_A,\gamma_A}))$.

Now, under the assumption (\ref{twoscale}) of local periodicity, one possible manner to get explicit expressions for the homogenized problem is to perform a formal two-scale analysis with the following Ansatz
\begin{equation} \label{Ansatz}
u^{\epsilon}(x) = u_{0}\left(x,\frac{x}{\epsilon}\right) + \epsilon u_{1}\left(x,\frac{x}{\epsilon}\right) 
+ \epsilon^{2} u_{2}\left(x,\frac{x}{\epsilon}\right) \dots
\end{equation}
where, for any $x \in \Omega$, the functions $u_{i}(x,\cdot)$ are Y-periodic. The first two terms of the Ansatz (\ref{Ansatz}) are shown to co\"incide with the $H^1$ approximation (\ref{corrector}) for $u^{\epsilon}$ \cite{Bensoussan,Allaire-92}. 

Inserting the Ansatz (\ref{Ansatz}) into equation (\ref{pbeps}) gives the following explicit expressions for the objects previously defined by the abstract homogenization result:
\begin{description}
\item[$\bullet$] the function $u_{0} = u^{\star}(x)$ does not depend on the fast variable $\epsilon^{-1}x$ and is the $L^2$ approximation for $u^\epsilon$ given by the convergence result (\ref{L2approx}),
\item[$\bullet$] the gradient $\nabla_y u_{1}(x,\cdot)$ linearly depends on $\nabla_x u^\star(x)$,
$$u_{1}\left(x,\frac{x}{\epsilon}\right) = \sum_{i=1}^{n} \partial_i u^\star(x) w_{i}\left(x,\frac{x}{\epsilon}\right) + \tilde{u}_{1}(x),$$
where $(w_{i}(x,\cdot))_{1\leq i\leq n}$ are $n$ Y-periodic \emph{cell functions},
\item[$\bullet$] the $n$ cell functions $w_{i}(x,\cdot)$, parametrized by their macroscopic position $x \in \Omega$, are solutions to the following $n$ \emph{cell problems},
\begin{equation} \label{pbcell}
-{\rm div}_{y}( \bar{\bar{A}}(x,y) \cdot [e_{i} + \nabla_{y} w_{i}(x,y) ] ) = 0 , \forall y \in \mathbb{R}^{n} \ ,
\end{equation} 
and the correctors $z^{\epsilon}_{i}$ now read $z^{\epsilon}_{i} = x_i + \epsilon w_{i}(x,x/\epsilon) $,
\item[$\bullet$] the entries $\left( \bar{\bar{A}}^{\star}(x)_{i,j} \right)_{1\leq i,j\leq n}$ of the homogenized matrix $\bar{\bar{A}}^{\star}$ can be explicitly computed with the cell functions $w_{i}(x,\cdot)$,
\begin{equation} \label{Astar}
\bar{\bar{A}}^{\star}(x)_{i,j}= \int_{Y}  \bar{\bar{A}}(x,y) [e_{i} + \nabla_{y} w_{i}(x,y)  ] \cdot e_{j} \ dy,
\end{equation}
\item[$\bullet$] the $H^1$ approximation for $u^\epsilon$ is now tractable and writes
\begin{equation} \label{approximate_corrector}
u^\epsilon = u^\star + \epsilon \mathop{\sum}_{i=1}^{n} w_i \partial_i u^\star + r_\epsilon\ ,
\end{equation}
where, provided $u^\star \in W^{2,\infty}(\Omega)$, the correction error $r_\epsilon$ can be estimated to locally scale as $\epsilon$ (far enough from the boundary layer), and to globally scale as $\sqrt{\epsilon}$,
\begin{eqnarray} \label{periodic_homogenization_error}
\| r_\epsilon \|_{H^1_{\Gamma_D}(\omega)} \leq C_1 \epsilon \| u^\star \|_{W^{2,\infty}(\omega)}, 
\forall \omega \Subset \Omega, \\
\| r_\epsilon \|_{H^1_{\Gamma_D}(\Omega)} \leq C_2 \sqrt{\epsilon} \| u^\star \|_{W^{2,\infty}(\Omega)},
\end{eqnarray}
with constants $C_i$ depending only on $\Omega$.
\end{description}

To sum up, the local periodicity assumption (\ref{twoscale}) allows to completely determine the homogenized problem through explicit two-scale expressions. The derivation of the homogenized equation in the case of locally periodic coefficients serves as a basis for many numerical homogenization strategies.

\subsubsection{Numerical homogenization strategies}

Under local periodicity assumption (\ref{twoscale}), a two-scale explicit homogenization strategy for a sequence of linear scalar elliptic PDEs like (\ref{pbeps}) reads as follows in the frame of Finite-Element approximations for the scalar elliptic problems (\ref{pbcell}) and (\ref{pbstar}).
\begin{algo}[Two-scale homogenization strategy] 
\label{algo::twoscale}
To homogenize the sequence of PDEs (\ref{pbeps}):
\begin{enumerate} 
\item solve the parametrized cell problems (\ref{pbcell}) at each point $x \in \Omega$ where the value of $\bar{\bar{A}}^{\star}(x)$ is necessary to compute the FE matrix of the homogenized problem (\ref{pbstar}),
\item store the functions $w_i$ for the future computation of the $H^1$ approximation of $u^{\epsilon}$,
\item assemble the FE matrix associated with the homogenized operator $- {\rm div}( \bar{\bar{A}}^{\star}\nabla \cdot)$,
\item solve the macroscopic homogenized problem (\ref{pbstar}),
\item build the $H^1$ approximation (\ref{approximate_corrector}) for $u^{\epsilon}$ with $u^{\star}$ and $w_i$.
\end{enumerate}
\end{algo}

On the one hand, in many practical situations, it is very common to assume that the tensors $\bar{\bar{A}}^{\epsilon}$ satisfy assumption (\ref{twoscale}). Indeed, in practice, $\bar{\bar{A}}^{\epsilon}$ is often known for some given $\epsilon=\epsilon_0$ only. So, the asymptotic structure $\bar{\bar{A}}^{\epsilon}$ of the problem with oscillating coefficients in $\Omega$ has to be constructed from the only member $\bar{\bar{A}}^{\epsilon_0}$. Now, to take advantage of the exact explicit expressions given by the two-scale analysis, it is preferable to build a family $\bar{\bar{A}}^{\epsilon}$ that satisfies assumption (\ref{twoscale}), when possible. Assumption (\ref{twoscale}) then seems fully justified for many applications from the practitioner's point of view. Then, the main numerical difficulty of the two-scale homogenization is the first step, that is the accurate computation of a large number of cell functions. This is the main issue addressed in this article.

On the other hand, for some applications where the heterogeneities are highly non-periodic, one may want to build the sequence $\bar{\bar{A}}^{\epsilon}$ differently, or even skip the explicit construction of the sequence $\bar{\bar{A}}^{\epsilon}$. For example, the actual construction process of heterogeneities may suggest another sequence $\bar{\bar{A}}^{\epsilon}$ for which the error estimation of the $H^1$ approximation would then be more precise and meaningful. Or it may seem too difficult to explicitly build such a sequence $\bar{\bar{A}}^{\epsilon}$ that satisfies (\ref{twoscale}) from the knowledge of some $\bar{\bar{A}}^{\epsilon_0}$ only. In such cases, many numerical homogenization strategies have been developped to treat the numerical homogenization of oscillating coefficients that are not locally periodic. 

To our knowledge, most of the existing numerical homogenization strategies may be classified in one of the two following categories. They either rely on different space assumptions than local periodicity for the oscillating coefficients ({\it e.g.} reiterated homogenization \cite{Lions-00}, stochastic homogenization \cite{Bourgeat-03}, deformed periodic coefficients \cite{Briane-94}, stochastically deformed periodic coefficients \cite{Blanc-06}), and still allow to derive exact (but not always fully explicit) expressions for the homogenized equation and the error estimate of the approximation. 

Or, the numerical homogenization strategies are much coarser and only rely on the assumption that explicit scale separation allow for the behaviour of the oscillating coefficients at some small $\epsilon$ to be numerically homogenized. Those strategies are then approximate-in-nature. They manage to approximate quite a large class of heterogeneous problems, but may be computationally very demanding. They may also lack sharp error estimates. Example are the Multiscale finite-element method (MsFEM) \cite{Hou-97,Allaire-04}, the Heterogeneous multiscale method (HMM) \cite{E-07}, or the recent variational approach for non-linear monotone elliptic operators proposed in \cite{Gloria-06-1}...

Now, in any of the two previously described situations where the numerical homogenization strategies require the computation of a large number of parametrized cell problems, the RB approach proposed here-after is likely to bring some additional computational efficiency. As a matter of fact, most numerical approximate homogenization strategies are only slight modifications of the exact two-scale homogenization strategy proposed above in the frame of local periodicity assumption (\ref{twoscale}), and they do require the computation of a large number of parametrized cell functions. For many mechanical applications, this owes to the assumed existence of a Representative volume elements (RVE), which leads to general cell problems at each point $x$ of the macro domain \cite{Gloria-06-2}. One simple example of a possibly approximate numerical homogenization strategy that requires the computation of a large number of parametrized cell functions is based on the following theorem, proved by Jikov {\it et al.} in \cite{Jikov-94}. 

\begin{theorem}
Let $\bar{\bar{A}}^{\epsilon}$ be a sequence of matrices in $L^{\infty}(\Omega,{\cal M}_{\alpha,\eta})$ that defines a sequence of linear scalar elliptic problems like (\ref{pbeps}). The H-limit of $\bar{\bar{A}}^{\epsilon}$ is the homogenized tensor $\bar{\bar{A}}^{\star}$. 

For any $x \in \Omega$ and $\epsilon > 0$, and any sufficiently small $h>0$, let us define a sequence of locally periodic matrices $\bar{\bar{A}}^{\epsilon}_h$ (in the sense of \eqref{twoscale}),
\begin{equation} \label{periodized}
\bar{\bar{A}}^{\epsilon}_h(x) = \bar{\bar{A}}^{\epsilon}(x+h[\epsilon^{-1}x])\ ,
\end{equation}
where $[\epsilon^{-1}x]$ denotes the integer part of $\epsilon^{-1}x$. 

Then, for every $1 \leq i \leq n$, there exists a unique sequence of periodic solutions $w_i^{\epsilon,h}(x,\cdot)$ in the quotiented Sobolev space $H^{1}_{\#}(Y)/ \mathbb{R}$ of $Y$-periodic functions in $H^{1}(Y)$ that satisfy the $n$ cell problems
\begin{equation} \label{generalizedcellpb}
 -{\rm div}( \bar{\bar{A}}^{\epsilon}(x+hy) \cdot [e_{i} + \nabla_{y} w_{i}^{\epsilon,h}(y) ] ) = 0, \forall y \in Y
\end{equation}
in the $n$-torus $Y =[0,1]^n$. 

For each point $x \in \Omega$ and $\epsilon > 0$, we define a matrix $\bar{\bar{A}}^{\star}_{\epsilon,h}$ made of elements
$$ \bar{\bar{A}}^{\star}_{\epsilon,h}(x)e_i\cdot e_j = \int_{Y} \bar{\bar{A}}^{\epsilon}(x+hy)  \cdot [e_{i} + \nabla_{y} w_{i}^{\epsilon,h}(x,y)  ] \cdot [e_{j} + \nabla_{y} w_{j}^{\epsilon,h}(x,y)  ] \ dy $$
for any $(i,j)$ in $\{1,2,\dots,n\}^2$. Then, there exists a subsequence $h'\rightarrow 0$ such that
$$\mathop{\rm lim}_{h' \rightarrow 0} \mathop{\rm lim}_{\epsilon \rightarrow 0} \bar{\bar{A}}^{\star}_{\epsilon,h'}(x) = \bar{\bar{A}}^{\star}(x)\ .$$
\end{theorem}

This theorem shows that, for any family $\bar{\bar{A}}^{\epsilon}$, it is always possible to approximate the exact homogenized problem with the same explicit expressions than those obtained under the local periodicity assumption \eqref{twoscale}, after local periodization of $\bar{\bar{A}}^{\epsilon}$ like in \eqref{periodized}.

So, considering the landscape for the homogenization theory as described above, among alternatives way of improving the numerical homogenization strategies, we choose here to concentrate on speeding up the numerical treatment of a large number of parametrized cell functions, rather than, for example, refining the approximations leading to explicit cell problems for larger class of oscillating coefficients.

In the sequel, for the sake of simplicity, we assume that the sequence of tensors $\bar{\bar{A}}^{\epsilon}$ satisfies assumption (\ref{twoscale}) and apply the RB approach to the two-scale numerical homogenization strategy (Algorithm \ref{algo::twoscale}). Yet, the RB approach may apply as well with any numerical homogenization strategy that consists of first approximating $\bar{\bar{A}}^{\epsilon}$ by some sequence of tensors $\bar{\bar{A}}^{\epsilon}_h$ like in \eqref{periodized}, the latter leading to an explicit approximation for the homogenized problem after solving parametrized generalized cell problems like (\ref{generalizedcellpb}). Let us now concentrate on decreasing as much as possible the computational cost of solving (\ref{pbcell}) for many parameter values $x \in \Omega$.

\subsection{The reduced-basis method}
\label{RB_context}

Two critical observations allow to think that an output-oriented model order reduction technique like the RB method is likely to improve the repeated numerical treatment of parametrized cell problems (\ref{pbcell}). First, only \emph{outputs} of the cell functions are required to solve the homogenized problem, and an {\it a posteriori} estimation gives sharp error bounds for those outputs. 

Second, as extensively discussed above, the numerous parametrized cell problems arising from numerical homogenization strategies can be solved independently for each value of the parameter. Thus, a computational procedure based on an \emph{offline/online approach} should naturally allow for a reduction of the computation time in the limit of many queries. In particular, a large number of (and theoretically, an infinity of) parametrized cell problems occurs in the limit $\epsilon \rightarrow 0$ of the homogenization strategies, in order to compute the homogenized problem \eqref{pbstar} with non-periodic coefficients. And the number of homogenized problems to compute and solve can also be very large in practice, for instance in the frame of parameter estimation and optimization problems.

These two observations motivate an RB approach for the parametrized cell problems (\ref{pbcell}), which should significantly decrease the expense of computations in terms of CPU time for the homogenization problems where the offline stage is short compared to the online stage, or where the offline stage is not even an issue (like in real-time engineering problems for instance) \cite{Nguyen-05}. We are now going to introduce the basics of the reduced-basis method, well known to experts, who may want to directly proceed to the Section \ref{RB_approach}.

\subsubsection{The parametrized cell problem}

Let $\X$ be the quotiented space $H^{1}_{\#}(Y) / \mathbb{R}$ of $Y$-periodic functions that belong to the Sobolev space $H^{1}(Y)$. The Hilbert space $\X$ is imbued with the $\tilde{H^{1}}(Y)$-norm 
$$\| u \|_{\X} = \left( \int_{Y} \nabla u\cdot \nabla u \right)^{1/2}$$ 
induced by the inner product $\big( u , v \big)_{\X} = \int_{Y} \nabla u \cdot \nabla v$ for any $(u,v) \in \X\times\X$. In the dual space $\X'$ of $\X$, the dual norm is defined for any $g \in \X'$ by 
$$\left\| g \right\|_{\X'} = \sup_{v \in \X} \frac{g(v)}{\|v\|_{\X}}\ .$$

For any $x \in \Omega$, we define:
\begin{description}
\item[$\bullet$] a continuous and coercive bilinear form in $\X \times \X$ parametrized by $x \in\Omega$,
$$ a(u,v;x) = \int_{Y} \bar{\bar{A}}(x,y) \nabla u(y) \cdot \nabla v(y) dy, \ \forall (u,v) \in \X\times\X\ , $$
for which $\alpha_A$ and $\gamma_A^{-1}$ are respectively coercivity and continuity constants, 
\item[$\bullet$] and $n$ continuous linear forms in $\X$ also parametrized by $x \in\Omega$,
$$ f_{i}(v;x)= - \int_{Y} \bar{\bar{A}}(x,y) e_{i} \cdot \nabla v(y) dy, \ \forall v \in \X\ , 1\leq i\leq n .$$
\end{description}

Now, for any integer $i$, $1\leq i\leq n$, the i-th cell problem (\ref{pbcell}) for the cell functions $w_{i}(x,\cdot)$ rewrites in the following weak form: Find $w_{i}(x,\cdot) \in \X$ solution for
\begin{equation} \label{pbcellvar}
a(w_{i}(x,\cdot),v;x) = f_{i}(v;x)\ , \forall v \in \X\ ,
\end{equation}  
where $x \in \Omega$ plays the role of a parameter.

We set ${\cal M}_i=\{w_{i}(x,y),x \in \Omega\}$ the solution subspace of the i-th cell problem (\ref{pbcellvar}) induced by the variations of $x$ in $\Omega$, and 
$$ {\cal M}=\{w_{i}(x,y),x \in \Omega,1\leq i\leq N\}=\mathop{\bigcup}_{i=1}^{N} {\cal M}_i $$ 
the global solution subspace, that is the reunion of all solution subspaces for all cell problems.

\begin{remark} Note at this point that ${\cal M}_i$ and ${\cal M}$ should be seen as spaces induced by the family of coefficients $\left(\bar{\bar{A}}^\star(x,\cdot)\right)_{x\in \Omega}$, and not $x$. It is indeed always possible (and often useful) to use other explicit quantities than $x$ as parameters to map ${\cal M}_i$ and ${\cal M}$, provided that the variations of the parameters inside a given range of values induce the same family of coefficients and the same corresponding cell functions than $x \in \Omega$.
\end{remark}

Besides, for the sake of simplicity in the presentation of the RB method, the tensor $\bar{\bar{A}}(x,\cdot)$ will be assumed symmetric in the following. Hence, in the computation of the homogenized tensor $\bar{\bar{A}}^\star(x)$, the only interesting output for the $n$ solutions $w_{i}(x,\cdot)$ is given by a symmetric matrix $s$ of size $n\times n$ (somewhat similar to a compliance in the terminology of mechanics). The entries $(s_{ij})_{1\leq i,j \leq n}$ of the matrix $s$ are given by
\begin{equation} \label{compliance}
s_{ij}(x) = -f_{j}(w_i(x,\cdot);x) = \int_{Y}  \bar{\bar{A}}(x,y) \nabla w_i(x,y) \cdot e_{j} dy\ .
\end{equation}
But note that the RB approach still applies with non-symmetric tensors $\bar{\bar{A}}(x,\cdot)$ modulo slight modifications\footnote{When the tensor $\bar{\bar{A}}(x,\cdot)$ is not symmetric, a dual problem, adjoint to the problem \eqref{pbcellvar}, is introduced. The dual problem can be solved similarly to the ``primal" problem \eqref{pbcellvar} with an RB method, in a dual RB projection space. Last, the output should be rewritten like $s$ plus an additional term that accounts for the residual error due to the RB projection of the equation \eqref{pbcellvar}. This primal-dual approach is more extensively described in \cite{Nguyen-05} for instance.}.

The purpose of the RB method is to speed up the computation of a large number of solutions $w_{i}(x,\cdot) \in \X$ of (\ref{pbcellvar}) for many parameter values $x \in \Omega$ while controlling the approximation error for the output $s$. 

\subsubsection{Principle of the reduced-basis method}

The purpose of most order reduction techniques like the RB method is to solve the weak form of a PDE like (\ref{pbcellvar}) through a Galerkin projection method with a Hilbertian basis that is ``adapted" to the solution subspace $\cal M$. 

For instance, a Hilbertian basis that is adapted to the equations (\ref{pbcellvar}) when $x \in \Omega$ and $1\leq i \leq n$ is an orthonormal family $(\xi_{j})_{j \in \mathbb{N}}$ (orthonormal with respect to the ambient inner product $(\cdot,\cdot)_{\X}$) such that 
\begin{itemize}
\item the ambient solution space $\X \subset \overline{{\rm span}\{\xi_{j},j \in \mathbb{N}\}}$ is separable,
\item and, for a finite $N$-dimensional subspace $\X_{N}={\rm span}\{\xi_{j},1\leq j\leq N\}$ of $\X$, the Galerkin approximations $w_{iN}(x,\cdot) \in \X_{N}$ for $w_{i}(x,\cdot)$ that satisfy, for any $x$ in $\Omega$ and $1\leq i\leq N$,
\begin{equation} \label{galerkin}
a(w_{iN}(x,\cdot),v;x) = f_{i}(v;x)\ , \forall v \in {\X}_{N}\ ,
\end{equation}  
are ``sufficiently" close to $w_{i}(x,\cdot)$ for a given tolerable precision. 
\end{itemize}
But the previous definition is only vaguely stated until the ``tolerable precision" is mathematically defined.

One possible way of defining sufficient precision is to control the approximation error for $w_{i}(x,\cdot) \in \X$ with the natural norm $\|\cdot\|_{\X}$ of the cell problem. The reduced-basis method rather proposes to control the approximation error for some linear output like $s$, which is not very different in the present case where $\bar{\bar{A}}(x,\cdot)$ is symmetric \footnote{In general non-symmetric cases, the approximation error for linear outputs like $s$ can be expressed as a product of two approximation errors, one for the parametrized solutions $w_{i}(x,\cdot)$ and another one for some dual quantity that is solution for the problem dual to \eqref{pbcellvar}. But here, because of the symmetry and of the specific nature of the output, the so-called compliance in reference to mechanics, the dual problem is unnecessary and the approximation error for the output can be directly expressed as the square of the approximation error for the cell function.}, as it will be made clearer in section \ref{sec::bounds}. 

The RB method is based on the computation of a basis for a Galerkin projection space in $\X$ that is ``adapted" to the cell problem (\ref{pbcell}) in the sense of the minimization of the output approximation error. The approximation errors are made explicit through rigorous {\it a posteriori} estimates, which allow to {\it a posteriori} certify the efficiency of the model order reduction, that is, the convergence of the RB method when the size $N$ of the Galerkin projection space increases.

\subsubsection{Practice of the reduced-basis method}

The RB method computes a basis for $\X_{N}$ from an approximation ${\cal M}_{p}$ of $\cal M$,
$$ {\cal M}_{p} =\{w_{i}(x,\cdot),x_k \in \D,1\leq i\leq n\}\ ,$$ 
induced by a discrete sample $\D = \{x_k,1 \leq k \leq p\}$, $p \times n > N$, of parameter values distributed over the parameter space $\Omega$ ($\D \subset \Omega$). This is termed as the \emph{offline} stage. Such a model order reduction is efficient if the tolerable precision for the output approximation error is reached with a $N$-dimensional ``adapted basis" when $N$ is mcuh smaller than $\cal N$, where $\cal N$ is the number of degrees of freedom necessary for a generic numerical method, like the FE method, to reach the same precision. We call reduced basis such an ``adapted basis" $(\xi_{j})_{1\leq j\leq n\times N}$. 

Now, the RB treatment of (\ref{pbcellvar}) begins with the computation of a sample of cell functions that induces ${\cal M}_{p}$. The cell functions of the finite space ${\cal M}_{p}$ should then be approximated before the model order reduction is possible. An accurate and generic numerical method\footnote{Note that the time needed to compute a (possibly large) sample of $p$ accurate FE approximations can also be an issue, that can be dealt with by a pre-processing stage according to the parametrization, as it will be made clearer in section \ref{sec::construction}.}, with a large number $\cal N$ of degrees of freedom, like an FE method with a fine mesh for $Y$ for instance, should be used at the beginning of the offline stage to compute an approximation ${\cal M}^{\cal N}_{p}$ for ${\cal M}_{p}$. 

Then, the sample of solutions ${\cal M}_{p}$ from which the basis $(\xi_{j})_{1\leq j\leq n\times N}$ is built should be as representative as possible of $\cal M$. In the absence of information, $\D$ should be chosen arbitrarily. So, it is not only quite often impossible to choose {\it a priori} the ``right'' parameter sample $\D$ for $\X_N$ in order to minimize the output approximation error over every $N$-dimensional vector subspace of $\cal M$, but besides, $\D$ should not be too large so that the reduced-basis construction is fast compared to the online stage. Hence the necessity for a reliable {\it a posteriori} control of the RB approximation method, which allows to build fast a reduced basis $(\xi_j)_{1\leq j\leq n\times N}$, as it will be seen in section \ref{sec::construction}.

After building a reduced basis for the vector field $\X_{N}$, the Galerkin projection method is applied to the weak form for (\ref{pbcellvar}) at any $x \in \Omega$. That is, in this \emph{online} stage, the previous reduced basis is assumed to span a vector field $\X_N$ sufficiently close to the solution manifold $\cal M$ of the parametrized PDE so that we can compute fast a sufficiently accurate Galerkin approximation in $\X_N$ for the solution of the parametrized PDE at any parameter value $x \in \Omega$.

\section{Reduced-basis approach for the cell problem}
\label{RB_approach}

The RB approach for equation (\ref{pbcellvar}) needs to {\it a posteriori} estimate the approximation error for Galerkin solutions of the cell problem. In a second stage, this allows for an {\it a posteriori} estimation of the approximation error on the output $s$.

\subsection{Error bounds for the cell problem}
\label{sec::bounds}

The purpose of this section is to derive the error bound \eqref{w_error_bound} for the cell funtioncs, solutions of equation (\ref{pbcellvar}). This allows to {\it a posteriori} estimate the approximation error for Galerkin solutions of the cell problem, and their outputs through the error bound \eqref{outputbound}. To this end, let us introduce the linear operator $T^x : \X \rightarrow \X$ so that, for any $u \in \X$ and $x \in \Omega$,
$$ (T^xu,v)_{\X} = a(u,v;x), \ \forall v \in \X \ .$$ 
The existence of such an operator directly leans on the Riesz-Fr\'echet representation Theorem in the Hilbert space $\X$.

For any $1\leq i\leq N$ and $x \in \Omega$, the Galerkin approximation error 
\begin{equation} \label{galerkinerror}
\| w_{i}(x,\cdot) - w_{i N}(x,\cdot) \|_{\X}
\end{equation}
can be bounded starting from the following equality,
\begin{equation} \label{eq1}
a( w_{i}(x,\cdot) - w_{i N}(x,\cdot),v;x) = f_{i}(v;x) -  a(w_{i N}(x,y),v;x) , \forall v \in \X,
\end{equation}
which is easily obtained by substraction of (\ref{galerkin}) from (\ref{pbcellvar}).

Let us define the parametrized bilinear residual forms $g_i$ in $\X\times\X$ such that, for all parameter values $x \in \Omega$ and $1\leq i\leq n$,
$$ g_i(u,v;x) = a(u,v;x)-f_i(v;x), \forall (u,v) \in \X\times\X\ .$$
Then, equation (\ref{eq1}) with $v = w_{i}(x,\cdot) - w_{i N}(x,\cdot)$ allows to immediately derive the following estimates through the dual norm of the residual linear form for $w_{i}(x,\cdot)$ defined in $\X$
$$ v \rightarrow g_i(w_{i N}(x,\cdot),v;x)\ .$$

First, owing to the coercivity of the bilinear form $a$, we obtain the lower bound:
\begin{equation} \label{ineq1} 
\alpha_A \| w_{i}(x,\cdot) - w_{i N}(x,\cdot) \|_{\X} \leq \| g_i(w_{i N}(x,\cdot),v;x) \|_{\X'}\ ,
\end{equation}
for the Galerkin approximation error.

Second, in view of the continuity of the bilinear form $a$, we obtain the superior bound:
\begin{equation} \label{ineq2} 
\| g_i(w_{i N}(x,\cdot),v;x) \|_{\X'} \leq \gamma_A^{-1} \| w_{i}(x,\cdot) - w_{i N}(x,\cdot) \|_{\X}\ , 
\end{equation}
for the Galerkin approximation error.

Finally, note that it is possible to compute the dual norm of the linear form
$$ v \rightarrow g_i(w_{i N}(x,\cdot),v;x) = - a( w_{i}(x,\cdot) - w_{i N}(x,\cdot),v;x) $$
using the Riesz-Fr\'echet representant $T^x( w_{i}(x,\cdot) - w_{i N}(x,\cdot) )$ in the Hilbert space \nolinebreak[4] $\X$, and that one can obtain numerical approximations for $\alpha_A$ and $\gamma_A^{-1}$, either using the spectral properties of the matrices resulting from the Galerkin projection in large generic solution spaces during the offline stage, or using properties of the parametrization like in section \ref{sec::parametrization}. So, the Galerkin approximation error (\ref{galerkinerror}) can be {\it a posteriori} bounded using estimations (\ref{ineq1}) and (\ref{ineq2}).

For $x \in \Omega$ and $1\leq i\leq N$, we define {\it a posteriori} estimators $\Delta_{N}(w_{i}(x,\cdot))$ for the Galerkin approximation errors (\ref{galerkinerror}), using the previous superior bounds, by
\begin{equation} \label{w_error_bound}
\Delta_{N}(w_{i}(x,\cdot)) = \frac{\|a(w_{i}(x,\cdot)-w_{i N}(x,\cdot),\cdot;x)\|_{\X'}}{\alpha_A}. 
\end{equation}

The effectivities $\eta_{N}(w_{i}(x,\cdot))$ corresponding to the estimators $\Delta_{N}(w_{i}(x,\cdot))$,
\begin{equation}
\eta_{N}(w_{i}(x,\cdot)) = \frac{\Delta_{N}(w_{i}(x,\cdot))}{\| w_{i}(x,\cdot)-w_{i N}(x,\cdot) \|_{\X}},
\end{equation}
satisfy the following inequalities independently of $N$,
\begin{equation} 
1 \leq \eta_{N}(w_{i}(x,\cdot)) \leq \frac{\gamma_A^{-1}}{\alpha_A},
\end{equation}
which shows the stability of the error estimator $\Delta_{N}(w_{i}(x,\cdot))$.

Last, Galerkin approximations for the homogenized and the output matrix write
\begin{eqnarray}
\bar{\bar{A}}^{\star}_{N}(x)_{i,j}& = & \int_{Y}  \bar{\bar{A}}(x,y) [e_{i} + \nabla_{y} w_{iN}(x,y)  ] \cdot e_{j} \ dy, \\
s^N_{ij}(x) & = & \int_{Y}  \bar{\bar{A}}(x,y) \nabla w_{iN}(x,y) \cdot e_{j} dy\ .
\end{eqnarray} 
The {\it a posteriori} superior bound $\Delta_{N}(w_{i}(x,\cdot))$ for the Galerkin approximation error (\ref{galerkinerror}) will now allow us to derive a simple superior bound for output approximation errors. Indeed, we have for any $1\leq i,j \leq N$ and $x \in \Omega$,
\begin{eqnarray*}  
\mid s_{ij}(x) - s_{ij}^{N}(x) \mid  & = & \mid f_{j}(w_{i}(x,\cdot)-w_{i N}(x,\cdot);x) \mid \\
& = & \mid a(w_{j}(x,\cdot),w_{i}(x,\cdot)-w_{i N}(x,\cdot);x) \mid \\
& = & \mid a(w_{j}(x,\cdot)-w_{jN}(x,\cdot),w_{i}(x,\cdot)-w_{i N}(x,\cdot);x) \mid \\
&\leq & {\alpha_A} \Delta_{N}(w_{i}(x,\cdot)) \Delta_{N}(w_{j}(x,\cdot))
\end{eqnarray*}
since $w_{j}(x,\cdot)$ and $w_{i}(x,\cdot)$ are solutions for (\ref{pbcellvar}), and $w_{i N}(x,\cdot)$ is solution for (\ref{galerkin}).

We are finally in possession of an {\it a posteriori} superior bound $\Delta_{ij,N}^s(x)$ for Galerkin approximations of the output $s_{ij}(x)$,
\begin{equation} \label{outputbound}
\Delta_{ij,N}^{s}(x) = \frac{ \|a(w_{i}(x,\cdot)-w_{i N}(x,\cdot),\cdot;x)\|_{\X'} \|a(w_{j}(x,\cdot)-w_{j N}(x,\cdot),\cdot;x)\|_{\X'} }{\alpha_A} \ .
\end{equation}
Numerical approximations for $\Delta_{ij,N}^s(x)$ will allow us to build fast a reduced basis for cell problems (\ref{pbcellvar}). Note that $\Delta_{ij,N}^s(x)$ scales as the product $\Delta_{N}(w_{i}(x,\cdot)) \Delta_{N}(w_{j}(x,\cdot))$, hence the interest of model order reduction techniques for solutions $w_{i}(x,\cdot)$ without much loss of precision for output $s(x)$.

\begin{remark} \label{product}
Note that for the output error bounds to scale like the square of the error bound for the cell functions, it has been essential to have the following orthogonality property for any $x \in \Omega$,
$$ a(w_{i}(x,\cdot)-w_{i N}(x,\cdot),w_{j N}(x,\cdot);x) = 0\ , \forall 1\leq i,j \leq N.$$
That is why we have chosen to build only one RB projection space $\X_{N}$, spanned by all the parametrized cell functions $w_i(x_k,y)$ when $1 \leq i \leq n$ and $x_k \in\D$. Yet, note that without this choice, the same scaling can still be obtained with $n$ distinct RB projection spaces $(\X_{iN})_{1\leq i\leq n}$ for each of the $n$ solution subspaces ${\cal M}_i$, provided one slightly modifies the definition of the output. Namely, another output matrix $\sigma$ and its RB approximation $\sigma^N$ should then be defined, starting from $s$ and $s^N$, by adding a residual error. Their entries read, for $1\leq i,j\leq n$,
\begin{eqnarray}
\sigma_{ij}(x) = -f_j(w_i(x,\cdot);x) + g_i(w_i(x,\cdot),w_j(x,\cdot);x) \\
\sigma_{ij}^N(x) = -f_j(w_{iN}(x,\cdot);x) + g_i(w_{iN}(x,\cdot),w_{jN}(x,\cdot);x)\ ,
\end{eqnarray}
where $\sigma = s$ because the tensor $\bar{\bar{A}}(x,\cdot)$ is symmetric, and $\sigma^N = s^N$ only when the RB projection space is the same for all solution subspaces ${\cal M}_i$ (as above). Interestingly, the same additional residual term in the output $\sigma$ also arises when the tensor $\bar{\bar{A}}(x,\cdot)$ is not symmetric. It is then evaluated with \emph{dual} cell functions, solutions for a problem dual to the cell problems \eqref{pbcellvar} \cite{Nguyen-05}.
\end{remark}

\subsection{The reduced-basis construction}
\label{sec::construction}

Let us choose a sample $\D = \{ x_{k}, 1 \leq k \leq p \}$ of $p$ values for the parameter $x$ in $\Omega$. Unless some physical properties of the system guides the choice for $\D$, the parameter values $x_k$ should be $p$ realizations of a random variable uniformly distributed over $\Omega$. To accurately solve equation (\ref{pbcellvar}) for each parameter value $x_k$, we choose an FE method with a $\cal N$-dimensional FE vector space $\X_{\cal N}$. Typically, $\cal N$ is very large for the FE approximations to be accurate. The $(n \times p)$ FE approximations $\big(w_{i{\cal N}}(x_k,\cdot)\big)_{i,k}$ for $\big(w_{i}(x_{k},\cdot)\big)_{i,k}$ span an ${(n \times p)}$-dimensional vector space $\X_{n \times p} \subset \X$ that contains the approximation 
$${\cal M}^{\cal N}_p=\{w_{i{\cal N}}(x_k,y),x_k \in \D\}$$
of the solution subspace $\cal M$. 

\begin{remark} \label{rem::preprocessing}
Computing the $(n \times p)$ FE approximations $\big(w_{i{\cal N}}(x_k,\cdot)\big)_{i,k}$ can become a cumbersome preliminary task to the reduced basis construction when the corresponding FE matrices are difficult to assemble. That is why the RB method also includes some pre-processing to assemble fast those FE matrices. Such a pre-processing is very simple in the case where the parametrization of the oscillating coefficients is affine (this terminology will be made clearer in section \ref{sec::parametrization}). It might be more difficult in other more general cases. In the present work, we only treat the affine case. Some more elaborate results in non-affine cases, that are based on the extrapolation method introduced in \cite{Barrault-04,Grepl-06} for instance, will appear in \cite{Boyaval-PhD}.
\end{remark}

First, in the offline stage of the RB approach, we would like to build a $N$-dimensional RB projection subspace $\X_N \subset \X$ that also contains a very close approximation of ${\cal M}^{\cal N}_p$, thus of $\cal M$. $\X_N$ will be spanned by a \emph{reduced basis} $(\xi_{j})_{1\leq j\leq N}$ made of $N$ vectors of $\X_{n \times p}$, with $N < n\times p$. Moreover, $N \ll {\cal N}$ should be sufficiently small for the model order reduction to allow a significant gain of computation time.

Then, in the online stage, for any $x \in \Omega$, $w_{i}(x,y)$ is to be approximated, using the RB method, by some vector $w_{iN}(x,y)$ in the Galerkin projection space $\X_{N}$ of size $N$ that writes 
$$w_{i N}(x,y) = \mathop{\sum}_{j=1}^{N} w_{i N j}(x) \xi_{j}(y)\ .$$

The reduced basis $\left( \xi_j(y) \right)_{1\leq j \leq N}$ of $\X_{N}$ is built in order to best control the approximation error for outputs through the {\it a posteriori} error bounds derived above. This is performed in the offline stage as follows.

\begin{algo}[Offline algorithm] We build a reduced basis $\left( \xi_j(y) \right)_{1\leq j \leq N}$ from ${\rm Span}\{w_i(x_k,\cdot),x_k \in \D,1 \leq i \leq n\}$ as follows:
\begin{enumerate}
\item for some couple $\left(k^0(1),i^0(1)\right)$, $1\leq k^0(1) \leq p$ and $1 \leq i^0(1) \leq n$, compute the accurate FE approximation $w_{i^0(1){\cal N}}(x_{k^0(1)},y)$ for $w_{i^0(1)}(x_{k^0(1)},y)$, element of $ {\cal M}^{\cal N}_p = \{w_{i{\cal N}}(x_k,y),x_k \in \D,1 \leq i \leq n\}$,
\item set $j=1$, $\displaystyle \xi_1(y) = \frac{w_{i{\cal N}}(x_k,y)}{\| w_{i{\cal N}}(x_k,\cdot) \|_{\X}}$, 
\item while $j < N$,
\begin{enumerate} 
\item compute for every $x_k \in \D$ and $1 \leq i \leq n$ the $(n\times p)$ RB approximations $w_{ij}(x_k,y) \in {\X}_{j} = {\rm span}\{ \xi_k, 1\leq k\leq j\}$ for the $n$ cell problems (\ref{galerkin}),
\item for $\displaystyle \left(k^0(j+1),i^0(j+1)\right) = \mathop{\rm argmax}_{1\leq k\leq p,1 \leq i \leq n} \frac{\Delta_{j}(w_{ij}(x_k,\cdot))}{\| w_{ij}(x_k,\cdot) \|_{\X}}$, compute the accurate FE approximation $w_{i^0(j+1){\cal N}}(x_{k^0(j+1)},y)$ for $w_{i^0(j+1)}(x_{k^0(j+1)},y)$, element of $ {\cal M}^{\cal N}_p = \{w_{i{\cal N}}(x_k,y),x_k \in \D,1 \leq i \leq n\}$,
\item set $\displaystyle \xi_{j+1}(y) = \frac{ R_{j+1}(y) }{\| R_{j+1}(y) \|_{\X}}$
where $R_{j+1}$ is the remainder of the projection on the $j$-dimensional reduced basis,
$$ R_{j+1}(y) = w_{i^0(j+1)}(x_{k^0(j+1)},y)-\mathop{\sum}_{k=1}^{j} (w_{i^0(j+1)}(x_{k^0(j+1)},\cdot),\xi_k)_{\X} \xi_k(y) ,$$
\item do $j=j+1$.
\end{enumerate}
\end{enumerate}
\end{algo}

\subsection{Convergence of the reduced-basis method for the cell problem}
\label{sec::apriori}

The {\it a priori} convergence of Galerkin approximations for solutions of continuous and coercive elliptic equations like (\ref{pbcellvar}) is classical. It usually relies on the following lemma (see {\it e.g.} \cite{Strang-73} for a proof).
\begin{lemma}[C\'ea Lemma]\label{lem::cea} For any $1\leq i\leq n$, let $w_{i}$ be the solution of (\ref{pbcellvar}) and $w_{i{N}}$ its approximation in some $N$-dimensional Galerkin projection space $\X_{N} \subset \X$. Then we have, for any $x \in \Omega$,
$$ \| w_{i}(x,y) -  w_{iN}(x,y) \|_{\X} \leq \sqrt{\frac{\gamma_A^{-1}}{\alpha_A}} \inf_{w(y) \in \X_{N}} \| w_{i}(x,y) -  w(y) \|_{\X} .$$
\end{lemma}

To conclude that RB approximations like $w_{i{N}} \in \X_{N}$ {\it a priori} converge to $w_{i}\in\X$ when $N\rightarrow\infty$, it is then usual to use Lemma \ref{lem::cea} in order to {\it a priori} prove the convergence of the approximation method.
\begin{lemma} \label{lem::apriori} If there exists a dense separable subspace $\cal V$ of $\cal M$ and an application 
$ r_{N}:{\cal V}\rightarrow {\X}_{N}$ such that
\begin{equation} \label{projection}
\mathop{\rm lim}_{N \rightarrow \infty} \| v - r_{N}(v) \| = 0, \forall v \in {\cal V}\ ,
\end{equation}
then, by C\'ea Lemma, for any $x \in \Omega$ and $1\leq i\leq n$, RB approximations $w_{i{N}}(x,\cdot)$ converge to $w_{i}(x,\cdot)$ in the following sense
\begin{equation} \label{RBconvergence}
\mathop{\rm lim}_{N\rightarrow \infty} \| w_{i}(x,y) -  w_{iN}(x,y) \|_{\X}\ .
\end{equation}
That is, for all $\epsilon > 0$, there exists a positive integer $N(\epsilon)$ such that, $ \forall x \in \Omega$ and $1\leq i\leq N$, 
\begin{equation} \label{globalconvergence}
 \| w_{i}(x,\cdot) -  w_{iN}(x,\cdot) \|_{\X} \leq \epsilon, \forall N\geq N(\epsilon)\ .
\end{equation}
\end{lemma}

Let us then naturally choose $\cal V = M$, and $r_{N}$ as the projection operator from $\cal M$ to $\X_{N}$ for the inner product in $\X$. Unfortunately, the convergence assumed in (\ref{projection}) can only be shown insofar as we have information about $\D$, which amounts to knowing how the parameter values are selected to build $\X_{N}$ as $N$ increases. Such an assumption is unrealistic since, to choose the right parameter values $x_k$ for $\D$, one should already know $\cal M$ or some spectral representation of it \cite{Maday-02}. So the scope of Lemma \ref{lem::apriori} seems strongly limited, as any {\it a prioi} analysis of the RB method in general.

As a matter fact, the RB method is a \emph{practical} method of order reduction and can only be {\it a posteriori} shown to converge using reliable and computationally unexpensive error bounds that can be evaluated along the RB approximations.

Note last that, by definition, the Galerkin projection space $\X_{N}$ is built to converge to the manifold ${\cal M}_{\cal N}=\{w_{i{\cal N}}(x,y),x \in \Omega\}$ induced by the FE approximations $w_{i{\cal N}}(x,y)$ in the sense that, for some given parameter $x \in \Omega$ and $1\leq i\leq n$, there exists for all $\epsilon > 0$ a positive integer ${\cal N}_i(\epsilon,x)$ such that
$$ \| w_{i}(x,y) -  w_{i{\cal N}}(x,y) \|_{\X} \leq \epsilon, \forall {\cal N}\geq{\cal N}_i(\epsilon,x)\ .$$

So, let us assume that we have an error estimate for $\mathop{inf}_{w(y) \in \X_{n\times N}} {\| w_{i}(x,y) -  w(y) \|_{\X}}$ that is global in parameter space $\Omega$ like in Lemma \ref{lem::apriori}, in the limit $N \rightarrow \infty$. Even then, on account of the pointwise convergence of FE approximations in parameter space, the RB method can only converge in the following sense 
\begin{equation} \label{actualRBconvergence}
\mathop{\rm lim}_{N\rightarrow \infty} \mathop{\rm lim}_{{\cal N}\rightarrow \infty}\| w_{i}(x,y) -  w_{iN}(x,y) \|_{\X} = 0\ ,
\end{equation}
where $w_{iN}(x,\cdot)$ implicitly depends on $\cal N$ and where the limits for $N$ and $\cal N$ are not reversible. Yet, if the error estimate is also global in parameter space $\Omega$ with respect to the limit ${\cal N} \rightarrow \infty$, then the limit for $N$ and $\cal N$ be inverted.

\subsection{Error estimate for the asymptotic $H^{1}$ homogenized solution}

In the frame of the two-scale homogenization strategy, the asymptotic $H^{1}$ homogenized approximation for $u^{\epsilon}(x)$ in the limit $\epsilon \rightarrow 0$ is
$$u_{0}(x) + \epsilon\ u_{1}\left(x,\frac{x}{\epsilon}\right) = u^{\star}(x) 
+ \epsilon \sum_{1 \leq i \leq n} w_{i}\left(x,\frac{x}{\epsilon}\right) \partial_{i} u^{\star}(x)\ ,$$
which strongly converges to $u^{\epsilon}(x)$ in $H^1_{\Gamma_D}(\Omega)$ when $\epsilon \rightarrow 0$ if $u^\star \in W^{2,\infty}(\Omega)$. 

In this approximation, the homogenized solution $u^\star$ is the solution to the variational formulation (\ref{pbstarvar}) of the homogenized equation (\ref{pbstar})
\begin{equation} \label{pbstarvar}
\int_\Omega \bar{\bar{A}}^\star \nabla u^\star \cdot \nabla v = \int_\Omega f v + \int_{\Gamma_N} v, \forall v \in H^1_{\Gamma_D}(\Omega)\ .
\end{equation}

But in practice, one can only compute an RB approximation $\bar{\bar{A}}^{\star}_{N}$ for $\bar{\bar{A}}^\star$, namely with entries
$$\left( \bar{\bar{A}}^{\star}_{N}(x) \right)_{i,j} = \left( \int_{Y} \bar{\bar{A}}(x,y)dy \right)_{i,j} - s_{ij}^N(x)\ ,$$
which should be taken into account to estimate the approximation error for the asymptotic $H^{1}$ homogenized approximation. The following lemma \eqref{lem::limits} will show how the {\it a posteriori} control of the RB approximation allows to control the approximation error for the asymptotic $H^{1}$ homogenized approximation with an RB approach.

Let us then define an approximation for the asymptotic $H^{1}$ homogenized approximation,
$$ u^{\star}_N(x) + \epsilon \sum_{1 \leq i \leq n} w_{iN}\left(x,\frac{x}{\epsilon}\right) \partial_{i} u^{\star}_N(x)\ , $$
where $w_{iN}$ is the RB approximation for $w_{i}$ defined in previous sections, and $u^\star_N$ is an approximation for $u^\star$ that is solution in $W_{h_{hom}}$ for the discrete variational problem
\begin{equation} \label{pbstarvarRB}
\int_\Omega \bar{\bar{A}}_N^\star \nabla u_N^\star \cdot \nabla v = \int_\Omega f v + \int_{\Gamma_N} v, \forall v \in W_{h_{hom}}\ ,
\end{equation}
with $W_{h_{hom}} \subset H^1_{\Gamma_D}(\Omega)$ a discrete FE Galerkin projection space associated with a mesh of size $h_{hom}$ for $\Omega$. We have the following result.

\begin{lemma} \label{lem::limits}
Assume that $\Gamma_{D}$ is a measurable subset of $\partial \Omega$ with a positive $(n-1)-$\nolinebreak[3]dimensional measure (when $n>1$) so that a Poincar\'e inequality holds for elements of the Sobolev space $H^{1}_{\Gamma_{D}}(\Omega) = \{ v \in H^{1}(\Omega), v\mid_{\Gamma_{D}} = 0\}$.

If the approximations $w_{iN}(x,\cdot)$ converge to $w_{i}(x,\cdot)$ for all parameter values $x$ in $\Omega$ in the sense
\begin{equation} \label{conv_assum}
\mathop{\lim}_{N\rightarrow \infty} \mathop{\max}_{1\leq i\leq n} \left\{ \mathop{\rm esssup}_{x \in \Omega} \| w_{i}(x,y) -  w_{i N}(x,y) \|_{\X} \right\} = 0 \ ,
\end{equation}
then the asymptotic $L^2$ homogenized approximation $u^\star_N$ converges to $u^\star$, and so does the approximation for the asymptotic $H^1$ homogenized approximation of $u^\epsilon$. That is, we have the two results
$$ \mathop{\lim}_{N \longrightarrow \infty} \mathop{\lim}_{\epsilon \rightarrow 0 }\left\| u^{\star}(x)-u_{N}^{\star}(x) \right\|_{L^{2}(\Omega)} = 0 $$
and
$$ \mathop{\lim}_{N \longrightarrow \infty} \mathop{\lim}_{\epsilon \rightarrow 0 }\left\| u^{\star}(x)-u_{N}^{\star}(x) + \epsilon \sum_{1 \leq i \leq n} \left( w_{i}\left(x,\frac{x}{\epsilon}\right) \partial_{i} u^{\star}(x) - w_{i N}\left(x,\frac{x}{\epsilon}\right) \partial_{i} u_{N}^{\star}(x) \right) \right\|_{H^{1}(\Omega)} = 0 $$
where the two successive limits cannot be inverted.
\end{lemma}

\begin{remark}
As explained in section \ref{sec::apriori}, the assumption \eqref{conv_assum} can barely be satisfied {\it a priori}. But in practice, the error bounds derived in the {\it a posteriori} analysis of section \ref{sec::bounds} allow to check this assumption. The numerical results of Section \ref{sec::numeric} even show that the convergence with respect to $N$ in \eqref{conv_assum} is exponential.
\end{remark}

{\it Proof of Lemma \ref{lem::limits}.} To show this result, let us define two quantities, 
$$ E^{u^\star(x)}_{N} = u^{\star}(x)-u_{N}^{\star}(x)$$
and
$$ E^{\nabla u^\star(x)}_{N} = \nabla_x (u^{\star} - u_{N}^{\star}) (x) + 
\mathop{\sum}_{i=1}^{n} \left( \nabla_{y} w_{i}\left(x,\frac{x}{\epsilon}\right) \partial_{i} u^{\star}(x) - \nabla_{y} w_{i N}\left(x,\frac{x}{\epsilon}\right) \partial_{i} u_{N}^{\star}(x) \right)\ .$$

The approximation errors for the asymptotic $L^2$ and $H^{1}$ homogenized approximation of $u^{\epsilon}(x)$ now respectively write
$$ \left\| u^{\star}(x)-u_{N}^{\star}(x) \right\|_{L^{2}(\Omega)} 
= \left\|E^{u^\star(x)}_{N}\right\|_{L^{2}(\Omega)} $$
and
\begin{eqnarray*} 
\left\| u^{\star}(x)-u_{N}^{\star}(x) + \epsilon \sum_{1 \leq i \leq n} \left( w_{i}\left(x,\frac{x}{\epsilon}\right) \partial_{i} u^{\star}(x) - w_{i N}\left(x,\frac{x}{\epsilon}\right) \partial_{i} u_{N}^{\star}(x) \right) \right\|_{H^{1}(\Omega)} \\
= \sqrt{ \left\| E^{u^\star(x)}_{N} \right\|^2_{L^{2}(\Omega)} + \left\| E^{\nabla u^\star(x)}_{N} \right\|_{L^{2}(\Omega)}^{2}  + \underset{\epsilon \rightarrow 0}O(\epsilon) }\ .
\end{eqnarray*}

Thus, the proof consists of the two successive results
\begin{equation} \label{res1} 
\mathop{\lim}_{N \longrightarrow \infty} \mathop{\lim}_{\epsilon \rightarrow 0 } \left\|E^{u^\star}_{N}\right\|_{L^{2}(\Omega)}=0
\end{equation}
and 
\begin{equation} \label{res2}
\mathop{\lim}_{N \longrightarrow \infty} \mathop{\lim}_{\epsilon \rightarrow 0 } \left\|E^{\nabla u^\star}_{N}\right\|_{L^{2}(\Omega)} = 0 \ .
\end{equation}

First, let us begin with properties of the homogenized tensor. On account of definition (\ref{Astar}), $\bar{\bar{A}}^\star(x)$ is a positive definite and continuous matrix. 

Indeed, for any $x \in \Omega$, $\bar{\bar{A}}^\star(x)$ is positive definite
$$ 0 < \alpha_A u \cdot u \leq \alpha_A \left( u \cdot u + \int_Y \left| \mathop{\sum}_{i=1}^{n} u_i[\nabla_y w_i(x,y)] \right|^2 dy \right) \leq \bar{\bar{A}}^\star(x) u \cdot u, \forall u \in \mathbb{R}^n$$ 
since $w_i(x,\cdot)$ is periodic.

And there exists a positive constant $\gamma^\star(x)$ such that $\gamma^\star(x)$ is a continuity bound for $\bar{\bar{A}}^\star(x)$
$$ \bar{\bar{A}}^\star(x) u \cdot u \leq \gamma_A^{-1} \left( u \cdot u + \int_Y \left| \mathop{\sum}_{i=1}^{n} u_i[\nabla_y w_i(x,y)] \right|^2 \right) \leq \gamma^\star(x) u \cdot u, \forall u \in \mathbb{R}^n$$
since the bilinear form in $ \mathbb{R}^n \times \mathbb{R}^n $
$$ (u,v) \rightarrow \int_Y \left(\mathop{\sum}_{i=1}^{n} u_i [\nabla_y w_i(x,y)]\right).\left(\mathop{\sum}_{j=1}^{n} v_j [\nabla_y w_j(x,y)] \right) dy$$
is clearly continuous. 

Moreover, we have uniform continuity because $\Omega$ is bounded. That is, there exists a real number $\gamma_{A^\star} > 0$ such that, for any $x$ in $\Omega$, $\gamma^\star(x) \leq \gamma_{A^\star}$.

Second, $u^\star \in W_{h_{hom}}$ and $u^\star_N \in W_{h_{hom}}$ satisfy variational formulations (\ref{pbstarvar}) and (\ref{pbstarvarRB}). We then have the following equality
\begin{equation*}
\int_\Omega \bar{\bar{A}}^\star \nabla u^\star \cdot \nabla v = \int_\Omega f v + \int_{\Gamma_N} v
= \int_\Omega \bar{\bar{A}}_N^\star \nabla u_N^\star \cdot \nabla v \ ,\forall v \in W_{h_{hom}}
\end{equation*}
that we rewrite with $v=(u^\star-u_N^\star)$
$$ \int_\Omega \bar{\bar{A}}^\star \nabla (u^\star-u_N^\star) \cdot \nabla (u^\star-u_N^\star)
= \int_\Omega (\bar{\bar{A}}_N^\star-\bar{\bar{A}}^\star) \nabla u_N^\star \cdot \nabla (u^\star-u_N^\star)\ .$$

Because of the coercivity of $\bar{\bar{A}}^\star(x)$, we finally have the inequality
$$ \alpha_A \|\nabla (u^\star-u_N^\star)\|_{L^2(\Omega)} \leq \|\bar{\bar{A}}_N^\star-\bar{\bar{A}}^\star\|_{\infty} \|\nabla u_N^\star\|_{L^2(\Omega)}\ .$$

Moreover, the Poincar\'e inequality for $u^{\star} - u^{\star}_{N}$ in 
$H^{1}_{\Gamma_{D}}(\Omega)$ writes as follows,
$$\|u^{\star}-u^{\star}_{N}\|_{L^{2}(\Omega)} \leq {\cal P} \|\nabla (u^\star-u_N^\star)\|_{L^2(\Omega)}, $$
with a certain constant $\cal P$ which only depends on $\Omega$. We have established an error estimate for $\|E^{u^\star}_{N}\|_{L^{2}(\Omega)}$.

Next, since $\bar{\bar{A}}^\epsilon(x)$ is a positive definite matrix for any $x$ in $\Omega$, we deduce the following inequality
$$ \alpha_{A} \|E^{\nabla u^\star}_{N}\|^2_{L^{2}(\Omega)} \leq 
\int_{\Omega} \bar{\bar{A}}\left(x,\frac{x}{\epsilon}\right) E^{\nabla u^\star(x)}_{N} \cdot E^{\nabla u^\star(x)}_{N} dx\ .$$

In the limit $\epsilon \rightarrow 0$, on account of the periodicity of $\bar{\bar{A}}(x,\cdot)$ and $\nabla_{y} w_{i N}(x,\cdot)$, the previous inequality rewrites
\begin{align*}
\mathop{\lim}_{\epsilon \rightarrow 0 } \left\|E^{\nabla u^\star}_{N}\right\|_{L^{2}(\Omega)}^{2} \\ 
\leq  \mathop{\iint}_{ \Omega \times Y} \frac{\bar{\bar{A}}(x,y)}{\alpha_{A}} \Big[\mathop{\sum}_{i=1}^{n} (e_i+\nabla_{y} w_{i}(x,y)) \partial_{i} u^{\star}(x) - (e_i+\nabla_{y}w_{i N}(x,y)) \partial_{i} u_{N}^{\star}(x) \Big]^{2} dy dx \ . 
\end{align*}

Last, the definition (\ref{Astar}) of the homogenized tensor $\bar{\bar{A}}^{\star}$ allows to rewrite the expression
$$ \int_{Y} \bar{\bar{A}}(x,y) \Big[\mathop{\sum}_{i=1}^{n} (e_i+\nabla_{y} w_{i}(x,y)) \partial_{i} u^{\star}(x) -  \nabla_{y} (e_i+w_{i N}(x,y)) \partial_{i} u_{N}^{\star}(x) \Big]^{2} dy $$
and we finally get the following error estimate
$$ \mathop{\lim}_{\epsilon \rightarrow 0 } \left\|E^{\nabla u^\star}_{N}\right\|_{L^{2}(\Omega)}^{2}
\leq \int_{\Omega} \frac{A^{\star}}{\alpha_{A}} \nabla(u^{\star} - u_{N}^{\star}) \cdot \nabla (u^{\star} - u_{N}^{\star}) 
+ \int_{\Omega} \frac{A^{\star} - A^{\star}_{N}}{\alpha_{A}} \nabla u_{N}^{\star} \cdot \nabla u_{N}^{\star}\ ,$$
the superior bound of which is itself superiorly bounded by
$$ \frac{\gamma_{A^\star}}{\alpha_{A}}  \| \nabla (u^{\star} - u^{\star}_{N}) \|_{L^{2}(\Omega)}^{2} 
+ \frac{1}{\alpha_{A}} \| \bar{\bar{A}}^{\star} - \bar{\bar{A}}^{\star}_{N} \|_{\infty} \| \nabla u^{\star}_{N} \|_{L^{2}(\Omega)}^{2}\ .$$

In the end, we have the following error estimates 
\begin{eqnarray} \label{uestimate}
\mathop{\lim}_{\epsilon \rightarrow 0 } \left\|E^{\nabla u^\star}_{N}\right\|_{L^{2}(\Omega)}^{2}
\leq \frac{1}{\alpha_{A}} \Big( \| \bar{\bar{A}}^{\star} - \bar{\bar{A}}^{\star}_{N} \|_{\infty} \frac{\gamma_{A^{\star}}}{\alpha_{A}} + 1 \Big) 
\| \bar{\bar{A}}^{\star} - \bar{\bar{A}}^{\star}_{N} \|_{\infty} \| \nabla u^{\star}_{N} \|_{L^{2}(\Omega)}^{2} \\ 
\mathop{\lim}_{\epsilon \rightarrow 0 } \left\|E^{u^\star}_{N}\right\|_{L^{2}(\Omega)}^{2} =
\left\|E^{u^\star}_{N}\right\|_{L^{2}(\Omega)}^{2} \leq \frac{\cal P}{\alpha_{A}} \| \bar{\bar{A}}^{\star} -\bar{\bar{A}}^{\star}_{N} \|_{\infty} \| \nabla u^{\star}_{N} \|_{L^{2}(\Omega)}^{2}\ .
\end{eqnarray}

They show that the asymptotic homogenized approximations converge if the approximate homogenized tensor $\bar{\bar{A}}^{\star}_{N}$ converges to $\bar{\bar{A}}^{\star}$.

Now, recall that the homogenized tensor $\bar{\bar{A}}^{\star}_{N}$ converges to $\bar{\bar{A}}^{\star}$ if the approximations $w_{iN}(x,\cdot)$ converge to the cell functions $w_{i}(x,\cdot)$ since we have already obtained the following error estimate
\begin{eqnarray} \nonumber
\| \bar{\bar{A}}^{\star} - \bar{\bar{A}}_{N}^{\star} \|_{[L^{\infty}(\Omega)]^{n\times n}} & = & \mathop{\max}_{1\leq i,j\leq n} \{ \mathop{\rm esssup}_{x \in \Omega} \mid s_{ij}(x) - s^N_{ij}(x) \mid \} \\
\nonumber
& \leq & \gamma_A^{-1} \mathop{\max}_{1\leq i\leq n} \left\{ \mathop{\rm esssup}_{x \in \Omega} \| w_{i}(x,y) -  w_{i N}(x,y) \|_{\X} \right\}^2
\label{Aestimate}
\end{eqnarray}
to derive error bounds for the output $s$. This concludes the proof of Lemma \ref{lem::limits}.$\ \Diamond$

\subsection{Practical influence of the parametrization}
\label{sec::parametrization}

This section is devoted to the pre-processing used by the RB method in order to fast assemble the FE and RB matrices corresponding to the projections of the variational formulation \eqref{pbcellvar} of the cell problem on the discrete FE and RB Galerkin approximation spaces.

Indeed, for a given family $\bar{\bar{A}}(x,y)$ of tensors and a given range $\Omega$ for parameter $x$ values, the solution subspace $\cal M$ for cell problems (\ref{pbcellvar}) is completely determined and fixed. Then, from the theoretical point of view, the way functions $w$ in $\cal M$ explicitly depend on some parameter $x \in \Omega$, which we call the parametrization, should not influence the efficiency of the RB method as a model order reduction technique, however it induces the solution subspace $\cal M$. But in practice, the explicit parametrization of $\cal M$ can significantly account for the efficiency of the RB method, because it greatly influences the practical assembling of the matrix and vectors in the Galerkin projection method.

Only piecewise-affine parametrizations (according to the terminology explained hereafter) are treated in this work, which allows for a fast, very accurate and simple pre-processing of the FE and RB matrices. But the RB method also adapts to other types of parametrizations (remember that we already refered to \cite{Boyaval-PhD} for more elaborate results in non-affine cases, based on the extrapolation method introduced in \cite{Barrault-04,Grepl-06} for instance).

In the case of an affine parametrization, the assembling of the matrix and vectors corresponding to the Galerkin projection of cell problems is always fast and easy. By affine parametrization of the cell problems, we mean that $\bar{\bar{A}}(x,\cdot)$ depends on the parametrization in an affine manner as follows: for any $x \in \Omega$,
\begin{equation} \label{affinedependence}
\bar{\bar{A}}(x,y) = \bar{\bar{A}}_0(y) + \mathop{\sum}_{q=1}^{Z} \Theta^q(x) \bar{\bar{A}}_q(y), \forall y \in Y
\end{equation}
where:
\begin{description}
\item[-] the matrix $\bar{\bar{A}}_0(y)$ defines a parameter-independent continuous and coercive bilinear form in $\X\times \X $,
$$ a_0(u,v) = \int_{Y} \bar{\bar{A}}_0(y) \nabla u(y) \cdot \nabla v(y) dy, \forall (u,v) \in \X \times \X\ ,$$
\item[-] the functions $\Theta^q:\Omega \rightarrow \mathbb{R}$ are parameter-dependent coefficient functions and
\item[-] the matrices $\bar{\bar{A}}_q(y)$ define parameter-independent continuous bilinear forms in $\X\times \X $,
$$ a_q(u,v) = \int_{Y} \bar{\bar{A}}_q(y) \nabla u(y) \cdot \nabla v(y) dy, \forall (u,v) \in \X \times \X\ . $$
\end{description}

With such affine parametrizations, the numerical RB treatment of cell problems is straightforward. Let us detail its implementation. 

First, we follow the offline algorithm presented in the section \ref{sec::construction}. At each step of the offline stage, a cell problem (\ref{pbcell}) for some parameter value in $\D = \{x_k,1\leq k\leq p\}$ is to be explicitly solved in order to build the reduced basis $\left(\xi_j(y)\right)_{1\leq j \leq N}$ for ${\rm span}\{w(x,y),x\in\D\}$. We choose to use an FE method for this initial step of the offline stage that consists of accurately solving the variational formulation (\ref{pbcellvar}) with conforming $\mathbb{P}_1$ Lagrange finite elements and a fine mesh for $Y$. 

The solution space $\X$ is discretized into the vector space ${\X}_{\cal N}$ of continuous, piecewise linear functions. Let ${\cal T}_{Y}$ be a conformal mesh for the $n$-torus $Y=[0,1]^n$ made of ${\cal N}_t$ elements $(\Sigma_k)_{1\leq k\leq {\cal N}_t}$ of size $h_Y$. We write $\phi_{k}$ the FE basis functions associated with the ${\cal N}$ nodes $y_k$ in $Y$. Now, for $0\leq q\leq Z$, we define the FE matrices $M_{q} \in \mathbb{R}^{{\cal N} \times {\cal N}}$ with entries
$$ \left(M_q\right)_{ij} = a_k(\phi_{i},\phi_{j}) $$
for any $1 \leq i,j \leq {\cal N}$, and $n$ FE data vectors $F_{q,l} \in \mathbb{R}^{\cal N}$ ($1\leq l\leq n$) with entries, for any $1 \leq i \leq {\cal N}$,
$$ \left(F_{k,l}\right)_i = \int_{Y}  \bar{\bar{A_k}}(y) e_l \cdot \nabla_{y} \phi_{i}(y) dy\ .$$
Then, for any $x$ in $\D$, $1\leq i\leq n$, we compute the FE approximate solution 
$$ w_{i{\cal N}}(x,y) \in {\X}_{\cal N} = \mathop{\sum}_{k=1}^{\cal N} w_{i{\cal N}k}(x) \phi_k(y) $$
for the cell problem (\ref{pbcellvar}), that satisfies
\begin{equation} \label{pbcellvarFE}
\left( M_0+\mathop{\sum}_{q=1}^{Z}\Theta^q(x) M_q \right) w_{i{\cal N}}(x,y_l) = \left( F_{0,i}+\mathop{\sum}_{k=1}^{Z}\Theta^q(x) F_{q,i} \right)_l
\end{equation}
at each node $y_l$, $1\leq l\leq {\cal N}$. The FE problem \eqref{pbcellvarFE} can then be fast and very accurately assembled through linear combinations of the matrices $(M_q)_{0\leq q\leq Z}$ and vectors $(F_q,i)_{0\leq k\leq Z}$, $1\leq i\leq n$, that are to be kept into memory.

Then, in the online stage, we would like to treat fast cell problems (\ref{pbcellvar}) for any parameter value $x \in \Omega$. Let us project the FE matrices and vectors on the RB space ${\X}_{n\times N}$, 
$$ M^{RB}_{q} = \xi^{t} M_{q} \xi , 0\leq q\leq Z$$
and
$$ F^{RB}_{q} =  \xi^{t} F_{q}, 0\leq q\leq Z\ ,$$
with $\xi$ the ${\cal N} \times (n N)$ matrix with columns $\xi_k$, $1\leq k\leq n\times N$. The RB approximation $w_{iN}(x,y) = \mathop{\sum}_{k=1}^{N} w_{iNk}(x) \xi_k(y)$ for cell function $w(x,y)$ is solution of the linear system
\begin{equation} \label{pbcellRB}
\left( M^{RB}_0+\mathop{\sum}_{q=1}^{Z}\Theta^q(x) M^{RB}_q \right) w_{i{\cal N}}(x,y_l) = \left( F^{RB}_{0,i}+\mathop{\sum}_{k=1}^{Z}\Theta^q(x) F^{RB}_{q,i} \right)_l\ ,
\end{equation}
that is fast assembled through linear combinations in the present affine case. And, for $1\leq i,j\leq n$, the outputs are easily given by 
$$ s_{ij}(x) =  \mathop{\sum}_{l=1}^{n N} \left( F^{RB}_{0,i}+\mathop{\sum}_{q=1}^{Z}\Theta^q(x) F^{RB}_{q,i} \right)_l w_{iNl}(x) .$$
Moreover, an error bound for each RB approximation can also be derived fast in the online stage, following an offline-online strategy similar to that applied to the RB output prediction \cite{Nguyen-05}.

So, affine parametrization obviously allows for a fast assembling of the RB matrix and vectors in (\ref{pbcellRB}). It can be considered as an ideal frame for an efficient RB method, because the CPU time for the online solution of one cell problem (\ref{pbcellRB}) actually scales like the CPU time for solving a linear system of size $N$, and because the offline stage is actually very short in comparison with the online computation of a large number of cell functions. The possibility of a similar gain of computation time is not as obvious in the case of non-affine parametrizations, when \eqref{affinedependence} is not valid anymore. Then, the assembling of matrices and vectors needs more elaborate techniques \cite{Barrault-04,Grepl-06}. Yet, an efficient RB approach is still possible for quite a few situations, as shown in \cite{Boyaval-PhD}.

We finally treat the case of piecewise affine parametrizations that can be recasted into the class of affine parametrizations. The pre-processing that we propose in this case lies on the fact that, in practice, the RB approximations are numerical approximations for FE approximations, and not for the ``true" cell functions. Then, to apply the RB method to the parametrized FE approximations, it is possible to consider a ``global" parametrization made of a parameter for the oscillating coefficients and of another parameter for the FE method (for instance, the geometrical features of the mesh).

We deal with oscillating coefficients $\bar{\bar{A}}(x,\cdot)$ parametrized in a piecewise affine manner as follows:
\begin{description}
\item[-] for each $x \in \Omega$, the cell $Y$ can be partitioned into $d$ non-overlapping $Y_k(x)$ open subsets ($d \in \mathbb{N}^\star$ should be fixed) such that $Y \subset \mathop{\bigcup}_{k=1}^d \overline{Y_k}(x)$,
\item[-] there exists $d$ non-overlapping reference open subsets $Y^0_k$ such that $Y \subset \mathop{\bigcup}_{k=1}^d \overline{Y^0_k}$,
\item[-] for each $x \in \Omega$, there exists $d$ affine homeomorphisms, $ \Phi_k(x,\cdot) : Y_k^0 \rightarrow Y_k(x)$, $1\leq k\leq d$,
\item[-] and for every $1\leq k\leq d$, the family of functions $\left( \bar{\bar{A}}(x,\Phi(x,\cdot)) \right)_{x \in \Omega}$ restricted to $Y^0_k$ can be parametrized in an affine manner as defined in \eqref{affinedependence} by
\begin{equation} \label{partdef}
\bar{\bar{A}}(x,\Phi(x,y)) =  \bar{\bar{A}}_0(y) + \mathop{\sum}_{q=1}^{Z} \Theta^q(x) \bar{\bar{A}}_q(y), \forall y \in Y^0_k\ .
\end{equation} 
\end{description}

The function $\Phi(x,\cdot)$, defined almost everywhere in $Y$ by
$$\Phi(x,y) = \Phi_k(x,y), \forall y \in Y_k^0, 1\leq k\leq d, $$
maps a ``reference" cell onto the cell with parameter value $x$. After the mapping, the family of cell problems defined with these piecewise affine oscillating coefficients can then be treated as if the parametrization was affine like in \eqref{affinedependence}, provided one take into account the stretching of the domain at each parameter value $x$.

For this, we define $2(Z+1)$ tensors of rank $3$, $\left( \bar{\bar{\bar{M}}}_{k} \right)_{1\leq k\leq Z+1}$ and  $\left( \bar{\bar{\bar{F}}}_{k} \right)_{1\leq k\leq Z+1}$, by:
\begin{eqnarray}
\bar{\bar{\bar{M}}}_k = \mathop{\sum}_{l=1}^{{\cal N}_t} \mathop{\sum}_{i=1}^{\cal N} \mathop{\sum}_{j=1}^{\cal N} \left(
\int_{\Sigma_l} \bar{\bar{A_k}}(y) \nabla_{y} \phi_{i}(y) \cdot \nabla_{y} \phi_{j}(y) dy 
\right) e_l \otimes e_i \otimes e_j, \\
\bar{\bar{\bar{F}}}_{k} = \mathop{\sum}_{l=1}^{{\cal N}_t} \mathop{\sum}_{i=1}^{\cal N} \mathop{\sum}_{j=1}^{\cal N} \left(
\int_{\Sigma_l}  \bar{\bar{A_k}}(y) e_i \cdot \nabla_{y} \phi_{j}(y) dy
\right) e_l \otimes e_i \otimes e_j\ .
\end{eqnarray}
An accurate pre-processing in the piecewise affine cases is then possible that assembles fast $\mathbb{P}_1$-FE matrices, and corresponding RB matrices, by adding a mapping step to the linear combinations of the affine cases. Namely, with the family of vectors $\left(\bar{V}(x)\right)_{x \in \Omega}$ that accounts for the stretching of the mesh elements,
$$ \bar{V}(x) =  \mathop{\sum}_{l=1}^{{\cal N}_t} det\left( (\nabla_y \phi(x,y)) |_{\Sigma_l} \right), \forall x \in \Omega,$$
we easily get the FE matrix for any parameter value $x$ in $\Omega$ through the formula 
$$ \bar{V} \cdot \left( \bar{\bar{\bar{M}}}_0 + \sum_{q=1}^Z \Theta^q(x) \bar{\bar{\bar{M}}}_q \right), $$
and so on for the RB matrix. Note also that the reduced basis should then be orthonormalized at each parameter value $x$ in $\Omega$, because the inner product matrix also changes for each $x$.

\section{Numerical results}
\label{sec::numeric}

We now show numerical results for the RB approximation of a seemingly non-affine two-dimensional problem that is brought back to the affine setting after mapping of the cell $Y$. We do not show the MP-RB treatment of more general piecewise continuous parametrizations in this work, but elementary results for one dimensional problems can be found in\nolinebreak[4] \cite{Boyaval-PhD}. The two-dimensional problem is chosen here to show the efficiency of the RB method in a classical situation for the homogenization theory. To fix ideas, it consists of homogenizing the conductivity of a heterogeneous composite material in a domain $\Omega$, where a two-dimensional matrix is full of inclusions with varying positions and conductivity properties.

\subsection{Definition of the problem}

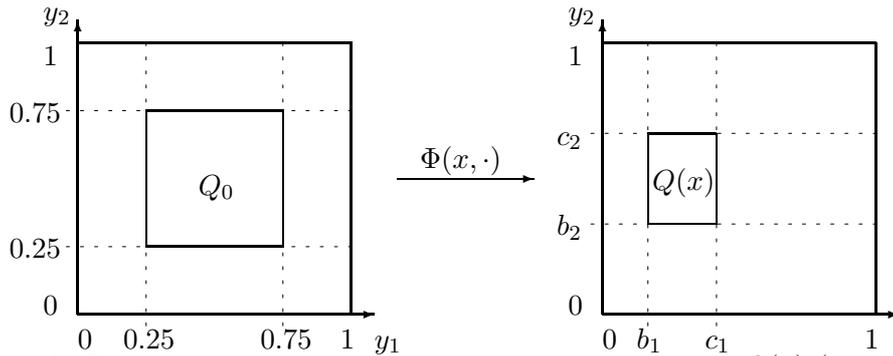
\begin{figure}[h!]
\centering
\setlength{\unitlength}{3mm}
\begin{picture}(36,13)

\linethickness{0.2mm}
\put(0,0){\vector(0,1){13}}
\put(0,-1.5){$0$}
\put(11.5,-1.5){$1$}
\put(13,-1.5){$y_1$}
\put(0,0){\vector(1,0){13}}
\put(-1.5,0){$0$}
\put(-1.5,11){$1$}
\put(-1.5,13){$y_2$}
\linethickness{0.25mm}
\multiput(0,0)(12,0){2}{\line(0,1){12}}
\multiput(0,0)(0,12){2}{\line(1,0){12}}

\linethickness{0.2mm}
\multiput(3,3)(6,0){2}{\line(0,1){6}}
\multiput(3,3)(0,6){2}{\line(1,0){6}}
\put(5.3,5.3){$Q_0$}

\dashline[+30]{.2}(3,-.5)(3,12)
\put(2,-1.5){$0.25$}
\dashline[+30]{.2}(9,-.5)(9,12)
\put(8,-1.5){$0.75$}
\dashline[+30]{.2}(-.5,3)(12,3)
\put(-3,2.5){$0.25$}
\dashline[+30]{.2}(-.5,9)(12,9)
\put(-3,8.5){$0.75$}

\linethickness{0.2mm}
\put(23,0){\vector(0,1){13}}
\put(23,-1.5){$0$}
\put(34.5,-1.5){$1$}
\put(36,-1.5){$y_1$}
\put(23,0){\vector(1,0){13}}
\put(21.5,0){$0$}
\put(21.5,11){$1$}
\put(21.5,13){$y_2$}
\linethickness{0.25mm}
\multiput(23,0)(12,0){2}{\line(0,1){12}}
\multiput(23,0)(0,12){2}{\line(1,0){12}}

\linethickness{0.2mm}
\multiput(25,4)(3,0){2}{\line(0,1){4}}
\multiput(25,4)(0,4){2}{\line(1,0){3}}
\put(25.2,5.5){$Q(x)$}

\dashline[+30]{.2}(25,-.5)(25,12)
\put(24.5,-1.5){$b_1$}
\dashline[+30]{.2}(28,-.5)(28,12)
\put(27.5,-1.5){$c_1$}
\dashline[+30]{.2}(22.5,4)(35,4)
\put(21,3.5){$b_2$}
\dashline[+30]{.2}(22.5,8)(35,8)
\put(21,7.5){$c_2$}

\linethickness{0.1mm}
\put(14,6){\vector(1,0){6}}
\put(15,6.5){$\Phi(x,\cdot)$}
\end{picture}
\caption[Geometry for 2D cells in numerical experiments]{\label{2Dcell} For each parameter value $x$, the cell with inclusion $Q(x)$ (on the right) is mapped through the piecewise affine homoemorphism $\Phi(x,\cdot)$ from a reference cell with inclusion $Q_0$ (on the left).}
\end{figure}

For $n=2$ and $f=0$, we supply the problem (\ref{pbeps}) with the mixed boundary conditions 
\begin{equation}
\text{(BC) } \left\{ \begin{array}{rcccl}
u^{\epsilon}(1,x_2) & = & 0 & = & u^{\epsilon}(x_1,1) \\
\bar{\bar{A}}^{\epsilon}\nabla u^{\epsilon}\cdot \bar{n}|_{(0,x_2)} & = & +1 & = &
\bar{\bar{A}}^{\epsilon}\nabla u^{\epsilon}\cdot \bar{n}|_{(x_1,0)}
\end{array} \right. \ .
\end{equation}
 
We define at each point $x \in \Omega$ a single rectangular inclusion $Q(x) \subset Y$ in the cell $Y = [0,1]^2$, $ Q(x) = \{ (y_1,y_2) | 0 < b_i(x) \leq y_i \leq c_i(x) < 1, 1\leq i\leq 2 \}$ (Fig.\ref{2Dcell}). We also write $\bar{\bar{I}}_2$ the second-order identity tensor and $\mathbf{1_{Q(x)}}$ the $Q(x)$-test function, such that, for every $y \in Y$, $\mathbf{1_{Q(x)}}(y)$ is one if $y \in Q(x)$ and zero otherwise. 

For all $x \in \Omega$, the oscillating coefficients $\bar{\bar{A}}^{\epsilon}(x) = \bar{\bar{I}}_2 + \bar{\bar{A}}_1(x,\epsilon^{-1} x)$ are locally periodic with a $Y$-periodic function $\bar{\bar{A}}_1(x,y) = \theta(x) \mathbf{1_{Q(x)}}(y) \bar{\bar{I}}_2$ that is constant inside and outside of the inclusion $Q(x)$. 

We want to homogenize the problem (\ref{pbeps}) associated with oscillating coefficients parametrized by the multiparameter $(b_1,c_1,b_2,c_2,\theta)(x)$, that is function of $x \in \Omega$ and takes value in $[.25-\delta;.25+\delta]^2 \times [.75-\delta;.75+\delta]^2 \times [-\theta^0;0]$, where $\delta \in ]0;.25[$ and $\theta^0 \in ]0;1[$. For the FE matrices to be easily assembled, we define a ``reference" cell problem with a centered inclusion $Q_0=[0.25;0.75]^2$ (Fig.\ref{2Dcell}). Then, at each point $x \in \Omega$, the inclusion $Q(x)$ can be mapped on $Q_0$ as explained in section \ref{sec::parametrization}.

\subsection{Offline computations}

\begin{figure}[h!]
\begin{tabular}{l l}
\includegraphics[width = 7cm]{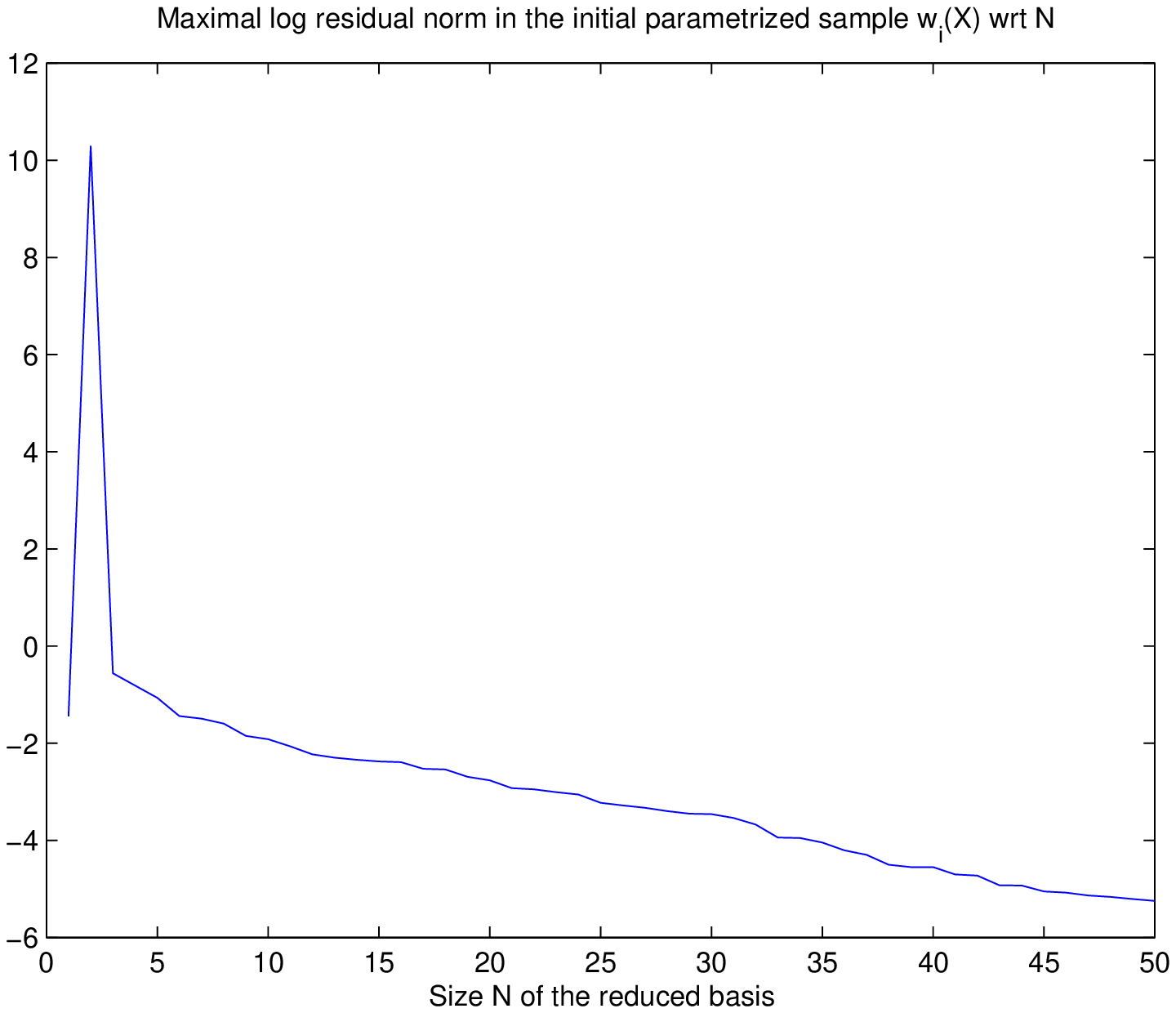} 
\includegraphics[width = 7cm]{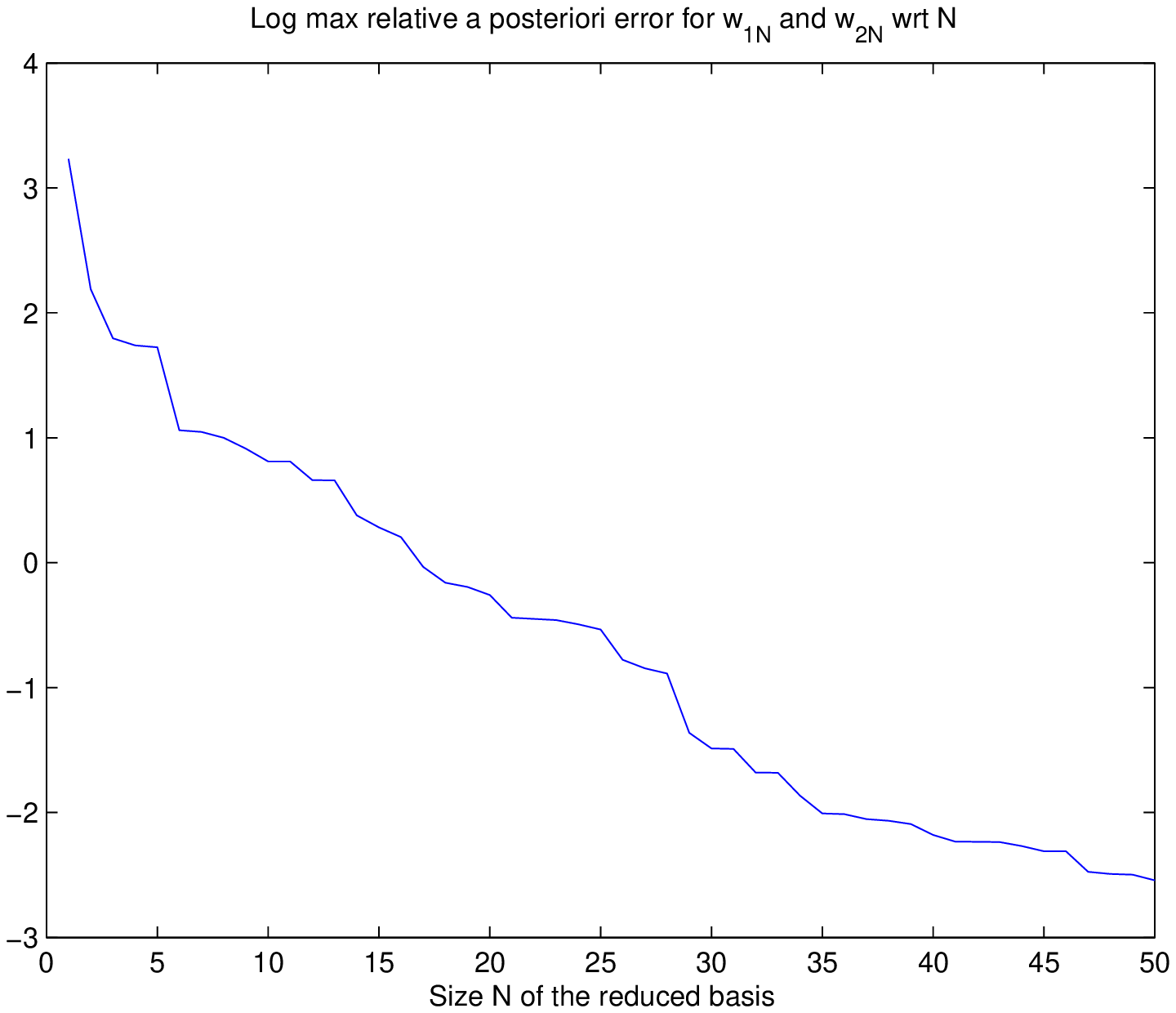} 
\end{tabular} 
\caption[RB efficiency after offline algorithm]{\label{2D_greedy} Maximal relative error bound $\displaystyle \mathop{\max}_{1\leq i\leq 2, x_k \in \D} \frac{\Delta_{N}(w_i(x_k,\cdot))}{\|w_i(x_k,\cdot) \|_{\X}}$ for the RB approximations $w_{iN}(x_k,\cdot)$ of the initial sample used for the RB construction (left picture), and for the RB approximations $w_{iN}(z_k,\cdot)$ of a test sample $z_k \in \Lambda$ (right picture), in log-scale with respect to the size $N$ of the growing reduced basis.}
\end{figure}

\begin{figure}[h!]
\begin{tabular}{l l}
\includegraphics[width = 7cm]{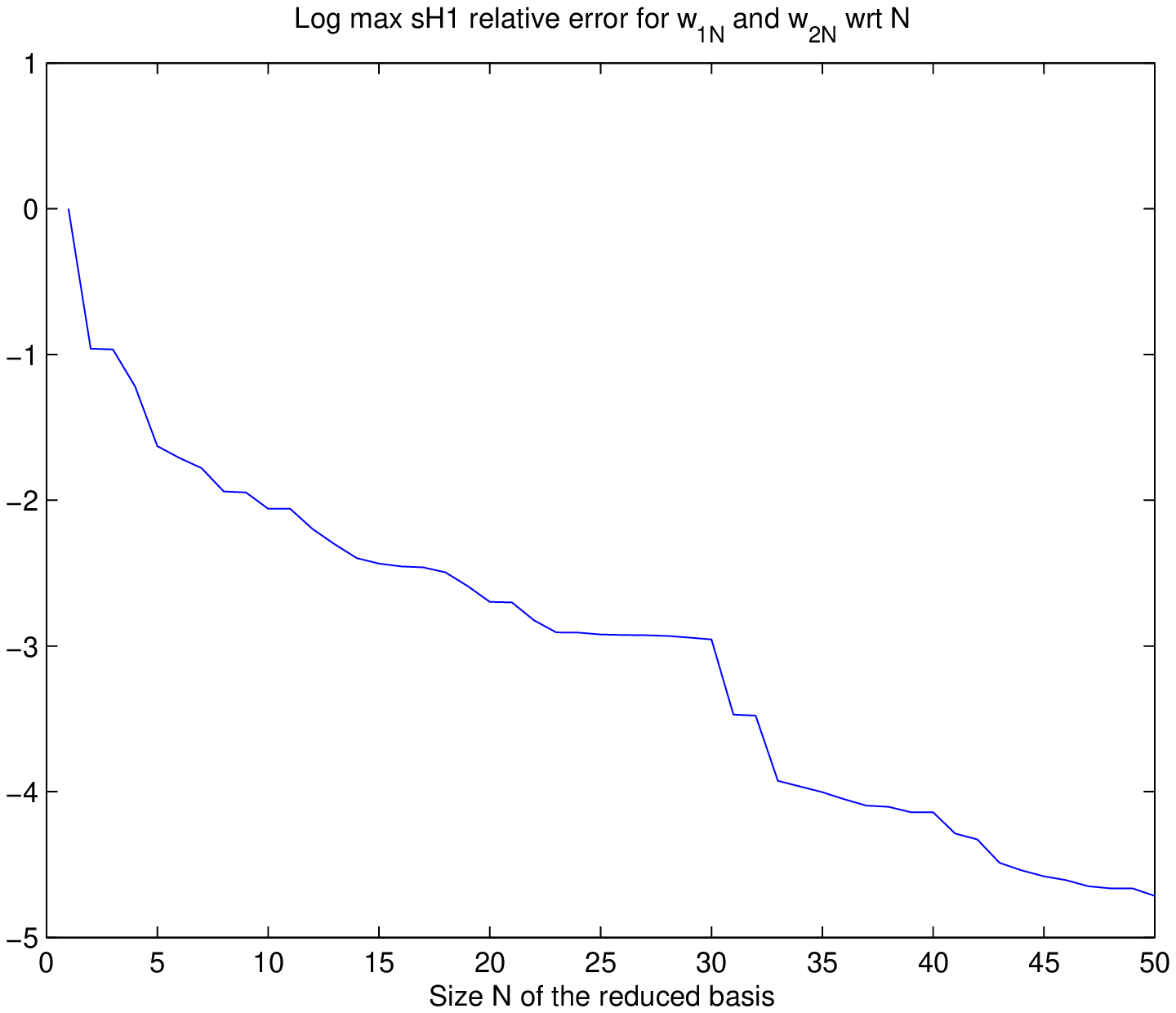} 
\includegraphics[width = 7cm]{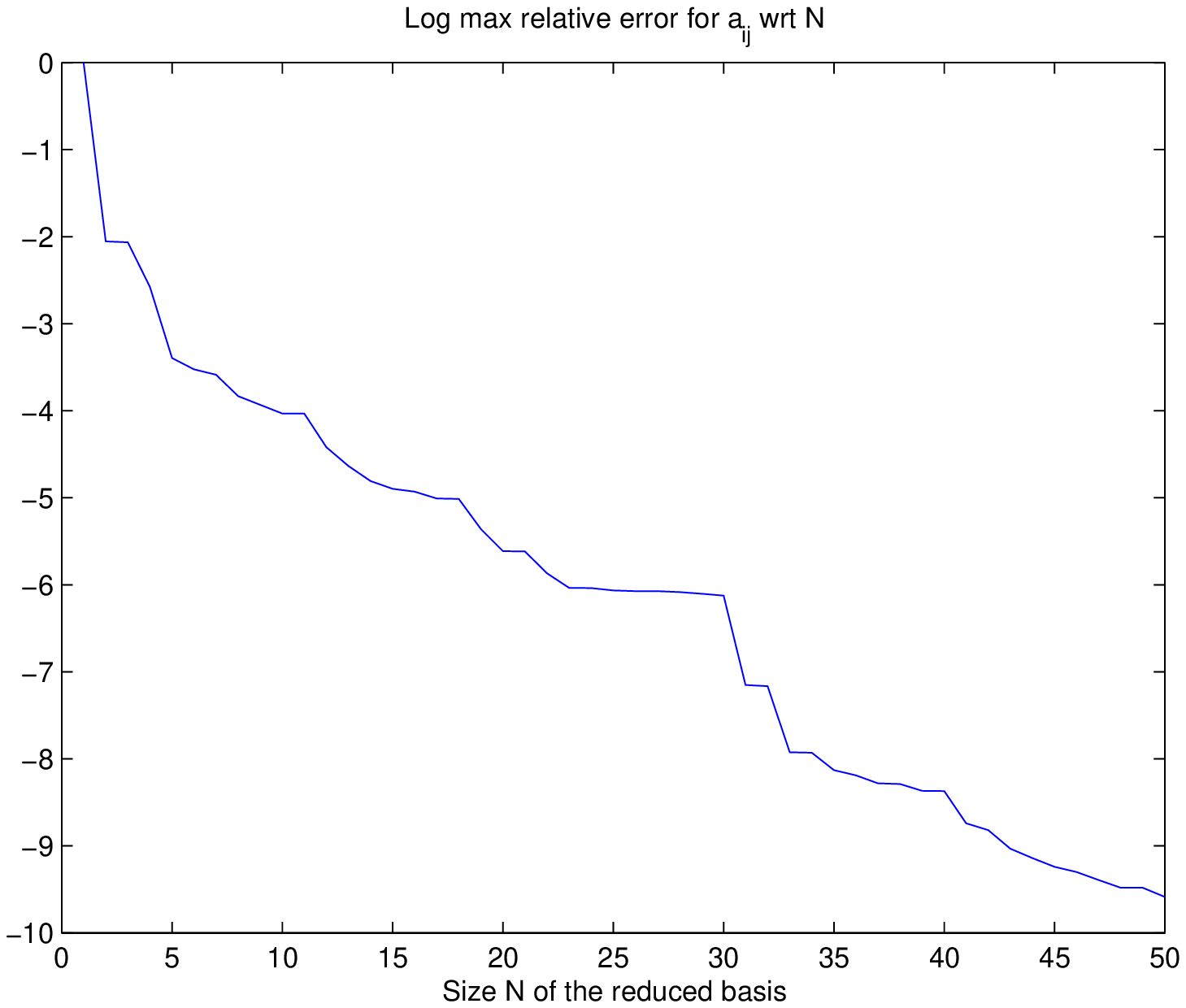}
\end{tabular} 
\caption[Test evaluation of the RB efficiency]{\label{2D_convergence} Maximal relative errors $\displaystyle \mathop{\max}_{1\leq i\leq 2, z_k \in \Lambda} \frac{\| w_i(z_k,\cdot)-w_{iN}(z_k,\cdot)\|_{\X}}{\| w_i(z_k,\cdot)\|_{\X}}$ (left picture) and $\displaystyle \mathop{\max}_{1\leq i,j\leq 2, z_k \in \Lambda} \frac{|s_{ij}(z_k)-s_{ij}^N(z_k)|}{|s_{ij}(z_k)|}$ (right picture) in log-scale with respect to $N$ for a random test sample $\Lambda$ of parameter values in $\Omega$.}
\end{figure}

A FE approach is developped for mapped cell problems in $Y$ with the ``reference" inclusion $Q_0$. More precisely, we use classical $\mathbb{P}_1$ simplicial Lagrange finite elements on a quadrangular, uniform and affine FE mesh, divided in isosceles triangles with base along direction $y_2=-y_1$ and size $h_Y$ in each direction $e_1$ and $e_2$. The mesh is fixed and adapted to the ``reference" domain in the sense that the boundaries of the inclusion $Q_0$ are multiples of $h_Y$.

We choose a random initial sample $\D$ of parameter values that is uniformly distributed over the multiparameter range. A reduced basis is then built for any parameter point $x$ after mapping with $\Phi(x,\cdot)$ the solutions $w_i(x_k,\Phi(x_k,y))$ selected by the offline algorithm of section \ref{sec::construction} for an initial sample $\D$ of $p=50$ parameter values. Numerical results are shown forh $\delta = .1$, $\theta^0 = .99$ and $h_Y = .1$ in Figures \ref{2D_greedy} and \ref{2D_convergence}. Note that the contrast between the coefficients inside and outside the inclusions can grow up to $1/100$, which is makes our expremient quite a stringent test.

The relative {\it a posteriori} error bounds for the RB approximations at the parameter values of the initial sample are computed at each step of the offline algorithm. The maximal error bound in this initial sample decreases exponentially with the size $N$ of the reduced basis (Fig. \ref{2D_greedy}). The effectivity of the {\it a posteriori} estimation is checked all along the RB construction (we found $\eta_N(w_i(x_k,\cdot)) \in ]1.4;3.5[$ for all $1\leq k\leq p$ and $1\leq i\leq n$ in the numerical experiment corresponding to Fig. \ref{2D_greedy}). Note that the offline algorithm selects (almost always alternatively) cell functions for the both cell problems, in direction $e_1$ and $e_2$. Then, for $N=2$, one cell function per direction only spans the reduced-basis, which strongly amplifies the RB approximation errors for the cell problem corresponding to the second direction represented in the reduced basis.

The reduced basis is then tested for another sample $\Lambda = \{z_k, 1\leq k\leq p\}$ of parameter values in $\Omega$. The maximal {\it a posteriori} error in this test sample still decreases very fast (exponentially with the size $N$ of the reduced basis), but the rate of decrease is slightly smaller than that of the initial sample used for the RB construction (Fig. \ref{2D_greedy}). This shows that the initial sample $\D$ was not an optimal choice to compute a reduced basis for any $x \in \Omega$, yet it still allows for efficient Galerkin approximations with any $\Lambda$. 

Besides, the {\it actual} output approximation error for this test sample scales as the square of the {\it actual} approximation error for cell functions (see Fig. \ref{2D_convergence}, obtained in the same numerical experiment than Fig. \ref{2D_greedy}). That is, the RB approximations are all the more efficient for the outputs, and the approximation errors scale like the error bounds derived in section \ref{sec::bounds}. The effectivities of the error bounds of section \ref{sec::bounds} are indeed hardly bigger than one (we found $\eta_N(w_i(z_k,\cdot)) \in ]1.3;3.9[$, for all $1\leq k\leq p$ and $1\leq i\leq n$, in the numerical experiment corresponding to Fig. \ref{2D_convergence})\footnote{Note that the maximal relative {\it a posteriori} error bound and the maximal {\it actual} error in Figures \ref{2D_greedy} and \ref{2D_convergence} are not obtained for the same parameter value $z_k$, hence the discrepancy between their ratio and the effectivity $\eta_N(w_i(z_k,\cdot))$ we measured.}.

\subsection{Online computations}

\begin{table}[h]
\centering
\begin{tabular}{|c|c|c|c|c|c|}
\hline
ratio & $(p=50)$ & Offline  & $( \frac{h_{hom}}{\epsilon} = \frac{3}{2} )$ & RB for $\bar{\bar{A}}_{N}^\star$ & FE for $\bar{\bar{A}}_{\cal N}^\star$ \\
$N / {\cal N}$ & $h_Y$ & algorithm & $\epsilon$ & (online) & (direct) \\
\hhline{|------|}
$1/5$ & $1 E^{-1}$ & 17 s  & $2 E^{-2}$ & 4+3 = 7 s & 27 s \\
\hhline{|------|}
$1/5$ & $1 E^{-1}$ & 15 s  & $2 E^{-3}$ & 410+330 = 740 s & 3100 s \\
\hhline{|------|}
$1/20$ & $5 E^{-2}$ & 42 s  & $2 E^{-2}$ & 16+10 = 26 s & 520 s \\
\hhline{|------|}
$1/20$ & $5 E^{-2}$ & 53 s  & $2 E^{-3}$ & 1600+1000 = 2600 s & 37000 s \\\hline\end{tabular}
\caption[CPU time for RB and FE approximate solutions]{\label{2D_time} CPU time (in seconds) needed by a Matlab code with an Intel Pentium IV processor (3.0 GHz/1 Go) to approximate the FE matrix for the homogenized problem either with a direct FE approach or with an RB method. In the RB approach, one has to take into account the RB construction (offline algorithm with a sample of $p$ parameter values), the online computation of one homogenized solution, plus possibly the online {\it a posteriori} estimation error (hence the two terms, solution+estimation, in the RB online column).}
\end{table}

\begin{table}[h]
$ \begin{array}{|c||c|c|c|c||}
\hline
( \frac{h_{hom}}{\epsilon} = \frac{3}{2} )
 & \|u^\epsilon-u^\star \|_{L^2}  & \|u^\star_{\cal N} - u^\star_{N}\|_{H^1} 
 & \|\nabla r_\epsilon\|_{L^2} & \|\nabla_y ( w_{i{\cal N}} - w_{iN}) \|_{L^2} \\
\hhline{|~|:====:|}
(\text{theory})
 & \leq C_1 \epsilon &  & \leq C_2 \sqrt{\epsilon} &  \\
\hhline{|~|----|}
h_Y = 1 E^{-1} & (\epsilon = 2.0 E^{-2}) & 1.2 E^{-4} & (\sqrt{\epsilon} = 1.4 E^{-1} ) & 2.9 E^{-2}  \\
\hhline{|~|----|}
h_Y = 1 E^{-1} & (\epsilon = 2.0 E^{-3}) & 4.7 E^{-3} & (\sqrt{\epsilon} = 4.5 E^{-2} ) & 1.0 E^{-2}  \\
\hhline{|~|----|}
h_Y = 5 E^{-2} & (\epsilon = 2.0 E^{-2}) & 3.1 E^{-3} & (\sqrt{\epsilon} = 1.4 E^{-1} ) & 8.6 E^{-5}  \\
\hhline{|~|----|}
h_Y = 5 E^{-2} & (\epsilon = 2.0 E^{-3}) & 1.1 E^{-3} & (\sqrt{\epsilon} = 4.5 E^{-2} ) & 3.0 E^{-2}  \\\hline
\end{array} $
\caption[RB approximation error]{\label{2D_approxerror} Theoretical correction error for the homogenized solution, and RB numerical approximation error for the homogenized solution when $\delta = .2$, $\theta^0 = .99$, $p=50$ and $N=20$.}
\end{table}

After building a reduced basis with the greedy algorithm from the previous FE approximations, we use the RB method to compute online RB approximations for cell functions as linear combinations of the RB basis functions. For this online stage, we develop an FE method for the homogenized problem \eqref{pbstar} and use classical $\mathbb{P}_1$ simplicial Lagrange finite elements on a quadrangular, uniform and affine FE mesh divided in isosceles triangles with base along direction $x_2=-x_1$ and size $h_{hom}$ in each direction $e_1$ and $e_2$.

The RB computations are performed in the step of the numerical homogenization strategy where the values of the homogenized coefficients are collected, as outputs of the cell functions, at some quadrature points in $\Omega$ that are necessary for the computations of the entries of the FE matrix in the homogenized problem \eqref{pbstar}. The CPU time needed for computing these outputs is compared between the RB and FE methods (Tab. \ref{2D_time}), where the RB method includes an online {\it a posteriori} estimation of its approximate solution. 

In calculations of Table \ref{2D_time}, the RB method has been applied with a reduced basis of size $N = 20$ starting from an initial parameter sample $\D$ of size $p = 50$. The main result of Table \ref{2D_time} is that the ratio of the RB computation time on the FE computation time scales like $N / {\cal N}$, the ratio of the numbers of degrees of liberty in the RB and (direct) FE methods.

So we can distinguish between two main regimes. The more unfavourable regime is the case of large ratios  $N / {\cal N}$, which corresponds to cases where one need only small precision for the correction \eqref{corrector} (large $h_Y$). Then, the RB method is likely to be faster than a direct FE method in the frame of many queries of the homogenized solutions. Note that in such a situation, the computation time spent by the offline algorithm is not even an issue. It is then possible to enlarge the initial parameter sample $\D$ (take a larger $p$). This increases the computation time spent by the offline algorithm in the RB construction but improves the quality of the RB approximations.

On the contrary, the favourable regime corresponds to small ratios $N / {\cal N}$, where the correction is sought very accurate (small $h_Y$). Then, the numerical results for the RB approximations of the cell functions show that there is an important gain of computation time, while there is no significant loss of numerical precision (Tab. \ref{2D_approxerror}). 

\section{Conclusion and perspectives}

We have shown in the present work that, for a prototypical class of parametrized cell problems (with piecewise affine oscillating coefficients), the reduced-basis approach applies and significantly reduces the time needed to compute a large number of parametrized cell problems in homogenization, in comparison with an FE method.

Some interesting questions concerning the extension of the RB approach in homogenization remain, mainly linked to the treatment of a larger class of parametrized cell problems: in particular, other geometries for more realistic cell problems should now be addressed, other boundary conditions for the cell problems (including the treatment of oversampling techniques), less regular oscillating coefficients (with many inclusions in varying amount). Also, the same questions as those examined in the present work for scalar elliptic equations could be asked for the Stokes-Darcy equations in porous media, or for the equations of linear elasticity in two- and three-dimensional contexts. Further developments of the RB methodology are then needed that may lead to interesting (fast) approaches in homogenization. Since a major issue in homogenization is the limitation of the time computation, speeding up the homogenization procedures could inevitably bring new possibilities of refinements (perhaps like reiterated oversamplings to improve the accuracy of the corrector term).

In any case, we believe that our result is interesting in the frame of many of the commonly used homogenization strategies, namely all those that ask for solving a computationally demanding number of parametrized cell problems. This is true provided the type of the parametrization can be handled with our RB approach. Among those homogenization strategies, the two-scale homogenization strategy is well known and much used in practice. That is why we have chosen this frame for our numerical experiments. But other homogenization strategies, which are used for non-locally-periodic oscillating coefficients, can also be treated with an RB approach.

For example, stochastic homogenization also asks for solving a large number of parametrized cell problems in the frame of local approximations of the homogenized tensor \cite{Bourgeat-03}. The homogenization of locally deformed oscillating coefficients, $ \bar{\bar{A}}^\epsilon(x) = \bar{\bar{A}}(\Phi^{-1}_x(\epsilon^{-1} x))$ with $\Phi$ a diffeomorphism, in the frame of deterministic homogenization \cite{Briane-94}
or of stochastic homogenization \cite{Blanc-06}, $ \bar{\bar{A}}^\epsilon(x) = \bar{\bar{A}}(\Phi^{-1}_x(\epsilon^{-1} x,\omega))$ with $\omega$ an element of a probability space, by nature, also demand for solving parametrized cell problems.

More general cases, often computationally demanding, also rely on the computation of a large number of cell problems, and offer a frame for an application of the RB approach. Among those homogenization strategies, the heterogeneous multiscale method (HMM), that averages over a large number of cell problems, could directly make use of our RB approach when cell problems are correctly parametrized. Another one, the multiscale finite-element method (MsFEM), also averages over numerous cell problems. Yet, the range of geometries for those cell problems is often larger, and it is still not obvious that model order reduction techniques may speed up the MsFEM computations.

Although we have not tested all the above mentioned possible improvement, we believe that our work is likely to improve a large number of existing homogenization strategies. Definite conclusions on the validity of our approach in such settings will hopefully be obtained soon.
 
\begin{mythanks}
I am indebted to C. Le Bris for suggesting to apply the RB methodology to homogenization problems, and to A.T. Patera for introducing me to the RB method. This work was initiated while I was visiting the Department of Mechanical Engineering at MIT (Boston, USA), and I would like to thank the group of A.T. Patera for their numerous advice. Last, this work has been possible thanks to grants from DARPA/AFOSR (FA9550-05-1-0114), from the Singapore-MIT Alliance and from CERMICS/ENPC (Marne-la-Vall\'ee, France).
\end{mythanks}

\bibliographystyle{siam}
\bibliography{./fig/refs_homogenization,./fig/refs_RB}

\end{document}